\documentclass[english,a4paper,11pt]{article}
\usepackage{amsmath,amssymb}
\usepackage{graphicx,url}
\usepackage{babel}

\newcommand{\RR}{\mathbb{R}}
\newcommand{\NN}{\mathbb{N}}

\newtheorem{Theorem}{Theorem}
\newtheorem{Proposition}{Proposition}

\begin{document}
\title{Elementary fractal geometry. \\ 
5. Weak separation is strong separation} 
\author{Christoph Bandt  \\ Mathematical Institute\\
University of Greifswald, Germany\\ \url{bandt@uni-greifswald.de} \and Michael F. Barnsley\\
Mathematical Sciences Institute\\
Australian National University, Canberra\\ \url{Michael.Barnsley@anu.edu.au}} 
\date{\today}
\maketitle
   
\begin{abstract}
For self-similar sets, there are two important separation properties: the open set condition and the weak separation condition introduced by Zerner, which may be replaced by the formally stronger finite type property of Ngai and Wang. We show that any finite type self-similar set can be represented as a graph-directed construction obeying the open set condition. The proof is based on a combinatorial algorithm which performed well in computer experiments.
\end{abstract}
 
\section{Overview}\label{intro}
A self-similar set is a nonempty compact subset $A$ of $\RR^d$ which is the union of shrunk copies of itself, as defined by Hutchinson's equation 
\begin{equation} A=\bigcup_{f\in F} f(A) \ .  
\label{hut}\end{equation}
Here $F$ denotes a finite set of contractive similarity mappings, called an iterated function system or IFS. The set $F$ fulfils the open set condition (OSC) if there exist an open set $U\subset \RR^d$ such that
\begin{equation} \bigcup_{f\in F} f(U) \subset U  \quad\mbox{ is a disjoint union of subsets of } U \ .   \label{osc}\end{equation}
This definition was introduced 1946 by Moran \cite{Mo} to show that the Hausdorff measure with respect to the similarity dimension $\alpha$ is positive on $A.$ It turned out to be the appropriate condition for many geometric studies of $A.$  The OSC says that the intersections $f(A)\cap g(A)$ with $f,g$ in $F,$ or in the generated semigroup 
\[ F^*=\bigcup_{k=0}^\infty F^k \quad \mbox{ with } \  F^k=\{ f_1f_2\cdots f_k\  |\ f_j\in F\} \]
are small: they are contained in the boundary of an open set. In algebraic terms, the OSC states that the IFS is discrete in the following sense \cite{BG}. There is a neighborhood $V$ of the identity map $id$ in the space of similitudes which does not contain maps of the form $f^{-1}g$ with $f,g\in F^*.$ For details see Chapter 4 in the recent book \cite{BSS23} and \cite{EFG3,EFG1,BP}.

Many examples of self-similar sets admit `exact overlaps' in the sense that $f(A)=g(A)$ for some $f,g\in F^*.$ Zerner \cite{Zer} noted that such coincidences of pieces, which obviously contradict the OSC, do not change the discrete character of the IFS. The weak separation condition (WSC) is satisfied if there is a neighborhood $V$ of $id$ which does not contain maps of the form $f^{-1} g$ with $f,g\in F^*,$ except for $id$ itself. 
In one dimension, Lau and Ngai \cite{LN99} studied this concept in a measure-theoretic setting. 

Zerner proved that self-similar sets with WSC have positive Hausdorff measure in some dimension $\beta$ which is smaller than the similarity dimension. While OSC needs to be checked, weak separation is satisfied for large classes of IFS; for example in the cases of crystallographic data \cite{Zer} and complex Pisot expansion factors \cite{EFG3}. 
Ngai and Wang \cite{NW} defined a finite type property which implies weak separation \cite{nguyen}, and gave an explicit formula for $\beta .$ All known examples with WSC, and a large class of one-dimensional examples fulfil the finite type condition \cite{Feng16,HHR21}, cf. \cite[Section 4.3]{BSS23}.  There are many recent studies on WSC attractors and their measures \cite{dajani21,DE11,DengWen20,FH15,HR22, KR16,Py21,TC15,Wu22}. \vspace{1ex}  

Here we show that weak separation is just a variant of the OSC. So in principle, there exists only one type of separation.  More precisely, we prove that
finite overlap type self-similar sets are graph-directed constructions \cite{MW}, abbreviated graph IFS or GIFS \cite{BV21,DE05}, with the open set condition.  The basic idea is to cut away the overlaps by going from an IFS to a GIFS.

To simplify the presentation, we assume that all similitudes $f\in F$ have the same contraction factor $r,$ that is $|f(x)-f(y)|=r|x-y|$ for all points $x,y\in\RR^d.$  We assume that the IFS has {\it finite overlap type:} there are only finitely many maps $h=f^{-1} g$ with $f,g\in F^k$ for some $k,$ such that $A\cap h(A)$ contains a whole piece $f'(A)$ for some $f'\in F^*.$  Compared with the literature \cite{BSS23, Feng16, HHR21, NW}, this is the weakest possible finite type condition. It is exactly the condition which Ngai and Wang needed when the gave their more complicated definition \cite[Section 2]{EFG3}. Finite overlap type implies WSC \cite{nguyen}, and all known examples with WSC have finite overlap type.  It is expected that the two properties coincide, which was proved for certain one-dimensional cases by Feng \cite{Feng16} and Hare, Hare and Rutar \cite{HHR21}. 

The GIFS in this paper have a simple form: a system of equations of the form \ref{hut} for $n$ attractors $B_k\subseteq A$ instead of the single attractor $A.$
\begin{equation} B_k=\bigcup_{f\in F_k} f(B_{j(f,k)})\quad \mbox{ for } k=1,...,n \ .  
\label{gifs}\end{equation}
Here $F_k\subseteq F$ so all equations involve only similitudes from $F.$ We can consider $j$ as a map $j:\bigcup_{k=1}^n F_k\times \{ k\} \to \{ 1,...,n\}$ which just means that each equation contains specific attractors $B_j$ for the maps $f\in F_k.$ This GIFS fulfils the OSC if there exist open sets $U_k, k=1,...,n$  such that
\begin{equation}\bigcup_{f\in F_k} f(U_{j(f,k)}) \subset U_k  \quad\mbox{ is a disjoint union of subsets of } U_k  \mbox{ for } k=1,...,n \ .   \label{gosc}\end{equation}

\begin{Theorem} \label{main} Let $F$ be an IFS of similitudes on $\RR^d,$ each with the same scaling factor, with the finite overlap type property. Let the attractor (self-similar set) be $A.$  Then $F$ can be extended to a GIFS of the form \eqref{gifs} with the OSC and $B_k\subseteq A.$
\end{Theorem}

This result goes beyond Hausdorff dimension, which can be easily determined from GIFS with the OSC \cite{MW}. The theorem shows that attractors with weak separation have the same structure as those with the OSC: they have a quite homogeneous modular structure consisting of pieces of positive Hausdorff measure in the respective dimension.

The GIFS will be explicitly constructed. The method was implemented on a computer. In the last part of the paper we handle new two-dimensional WSC examples of moderate complexity while the literature so far has been focussed on one-dimensional attractors \cite{BSS23, Feng16, HHR21, NW}. Treating self-similar sets of even larger complexity is a challenge concerning both the computer implementation and the mathematical problem of minimising the number of equations in a GIFS. 

Finite automata play the main part in this paper. They have been around in dynamical systems and fractal geometry for many years, without explicit mention. Finite type shifts and sofic shifts are automata-generated symbolic data \cite[Chapter 3]{LM95}. ``The concept of a `sophic system' is simply the ergodic theorist's specialized name for what is known in other branches of science as a finite state automaton'' \cite[Preface]{epstein}. A `graph-directed construction' \cite{MW} or GIFS is the addressing of fractal attractors by such automata-generated symbolic data. A `neighbor graph' \cite{EFG3,EFG4,BM09} or `overlap graph' addresses the dynamical boundary, or the overlaps, by automata-generated symbols.  In automata-theoretic terms, the construction of this paper describes a transducer from the overlap automaton of a WSC attractor to a GIFS automaton of OSC attractors. This means that our approach to WSC does not depend on similitudes. It can also be realized in an affine, conformal or topological setting, cf. \cite{EFG4,KR16,TC15}. 
Moreover, our computer experiments in Section \ref{tile} lead to questions of automata theory.

We start with two simple examples illustrating the basic idea. In Section \ref{ft} we introduce the overlap graph as our main tool. Section \ref{comp} compares this method with that of Ngai and Wang.  The GIFS construction is presented in Section \ref{const} and worked out for an example. The proof of the OSC for Theorem \ref{main} in Section \ref{prof} uses an idea of Schief \cite{Sch} and the concept of a central open set \cite{BHR}. Detailed examples in Sections \ref{deta} and \ref{deta2} and results of computer experiments in Sections \ref{prob} and \ref{tile} show that our method works well for fractals and tiles of moderate complexity. 

All fractal figures in this paper were produced by Mekhontsev's IFStile package \cite{M}. Overlap graphs were treated with MATLAB.  The first named author gratefully acknowledges an encouraging e-mail conversation with Dylan Thurston.  

\section{The basic idea} \label{basic}
\begin{figure}[h!t] 
\begin{center}
\includegraphics[width=0.3\textwidth]{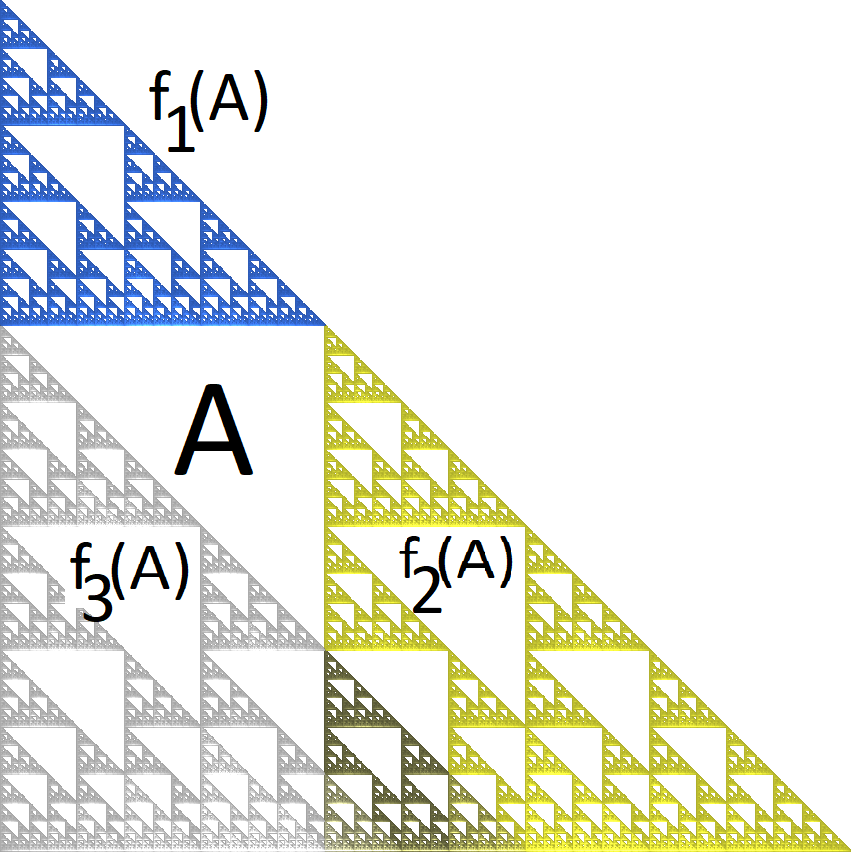}
\includegraphics[width=0.3\textwidth]{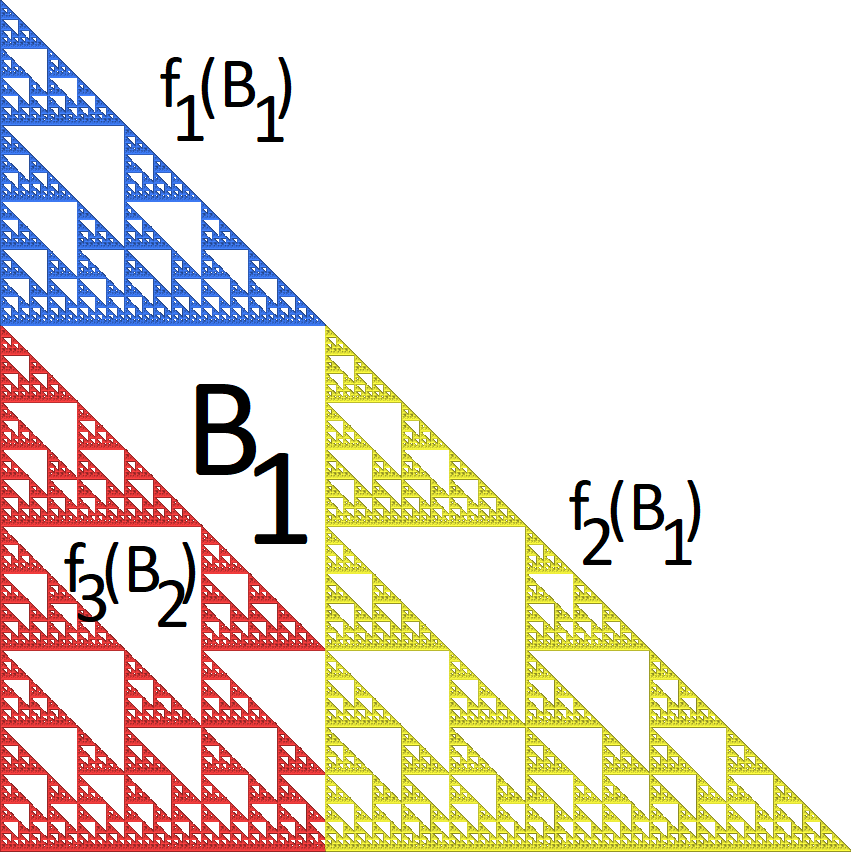} 
\includegraphics[width=0.19\textwidth]{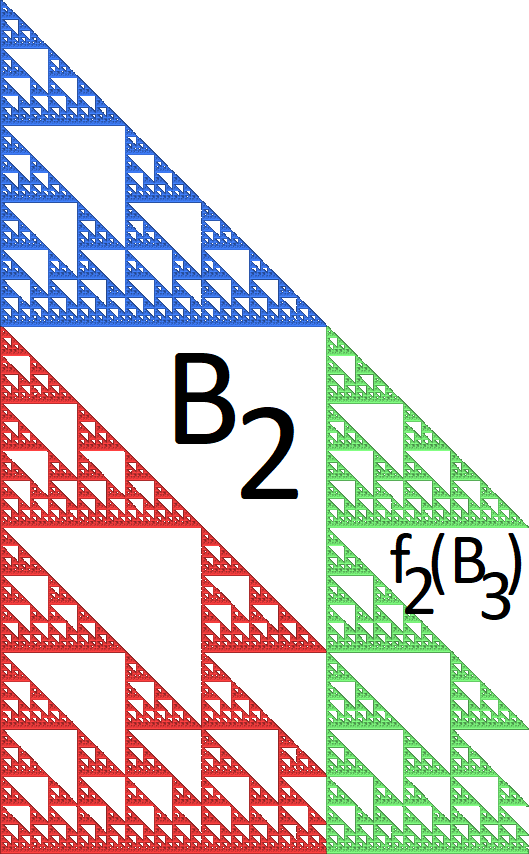} \
\includegraphics[width=0.118\textwidth]{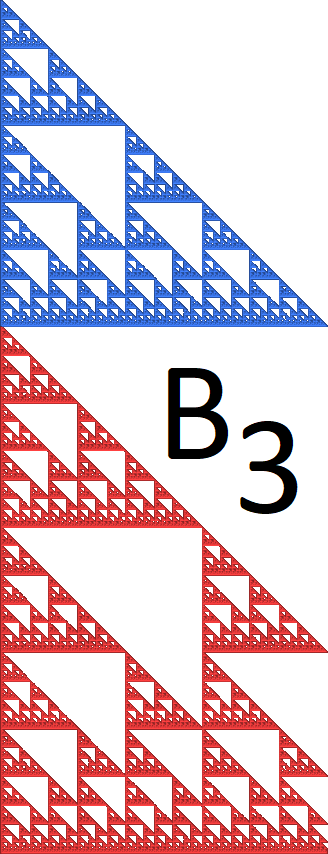}
\end{center}
\caption{The example of Ngai and Wang and an associated non-overlapping GIFS.}\label{top1}
\end{figure}

Figure \ref{top1} shows the only two-dimensional example of Ngai and Wang \cite[Example 5.4]{NW}, reproduced in \cite[Example 4.3.7]{BSS23}. Our approach works although not all maps have the same factor. The overlapping attractor $A=f_1(A)\cup f_2(A)\cup f_3(A)$ is generated by the IFS  
\begin{equation}
f_1(x,y)=(t^2x,t^2y+t),\ f_2(x,y)=(tx+t^2,ty), \ f_3(x,y)=(tx,ty) \ \mbox{ with } t=\textstyle\frac12 (\sqrt{5}-1)\ .
\label{nwifs}\end{equation}

We eliminate the overlap in the leftmost part of Figure \ref{top1} by introducing copies $B_k$ of $A$ with cut-away overlaps. Let $B_1=A$ with a non-overlapping representation. We replace the left part $f_3(B_1)$ of $B_1$ by $f_3(B_2)$ where $B_2$ is chosen so that $f_3(B_2)$ (drawn in red) and $f_2(B_1)$ touch each other within a vertical line segment.  
Actually, $B_2$ is the closure of $B_1\setminus f_3^{-1}f_2(B_1).$ Now $B_2$ is
a union of three pieces $f_1(B_1), f_3(B_2)$ and $f_2(B_3)$ where $B_3$ (drawn in green) is a new set smaller than $B_2.$ Finally, $B_3=B_2\cap f_2^{-1}(B_2)$ is the union of only two pieces, and we obtain the GIFS
\begin{equation}
B_1=f_1(B_1)\cup f_2(B_1)\cup f_3(B_2), \ B_2=f_1(B_1)\cup f_2(B_3)\cup f_3(B_2), \ B_3=f_1(B_1)\cup f_3(B_2)\, .
\label{nwgifs}\end{equation}
The IFS \eqref{nwifs} together with the equations \eqref{nwgifs} determine the GIFS. We did not calculate anything. We just looked at the structure of the overlap. This combinatorial approach will work in general.

The dimension $\beta$ of the $B_k$ can be obtained from \eqref{nwgifs} by the method of Mauldin and Williams \cite{MW}. For the present example there is a simpler trick. We know that the $\beta$-dimensional Hausdorff measure $\mu$ is positive and finite, and for a similitude $f$ with factor $t$ we have $\mu (f(E))=t^\beta \mu(E).$ We can assume $\mu(A)=1$ and write $s=t^\beta$ so that $\mu(f_2(A))=\mu(f_3(A))=s$ and $\mu(f_1(A))=s^2.$ The second and third piece overlap in a piece of third level so that $\mu(f_2(A)\cap f_3(A))=s^3.$ The measure of $A$ is the sum of the measure of the pieces $f_j(A)$ minus the measure of the overlap: $1=s^2+s+s-s^3.$
The root in $[0,1]$ is $s\approx 0.445,$ which implies $\beta =\frac{\log s}{\log t}\approx 1.682,$ simplifying the calculation in \cite{BSS23,NW}.

\begin{figure}[h!t] 
\begin{center}
\includegraphics[width=0.85\textwidth]{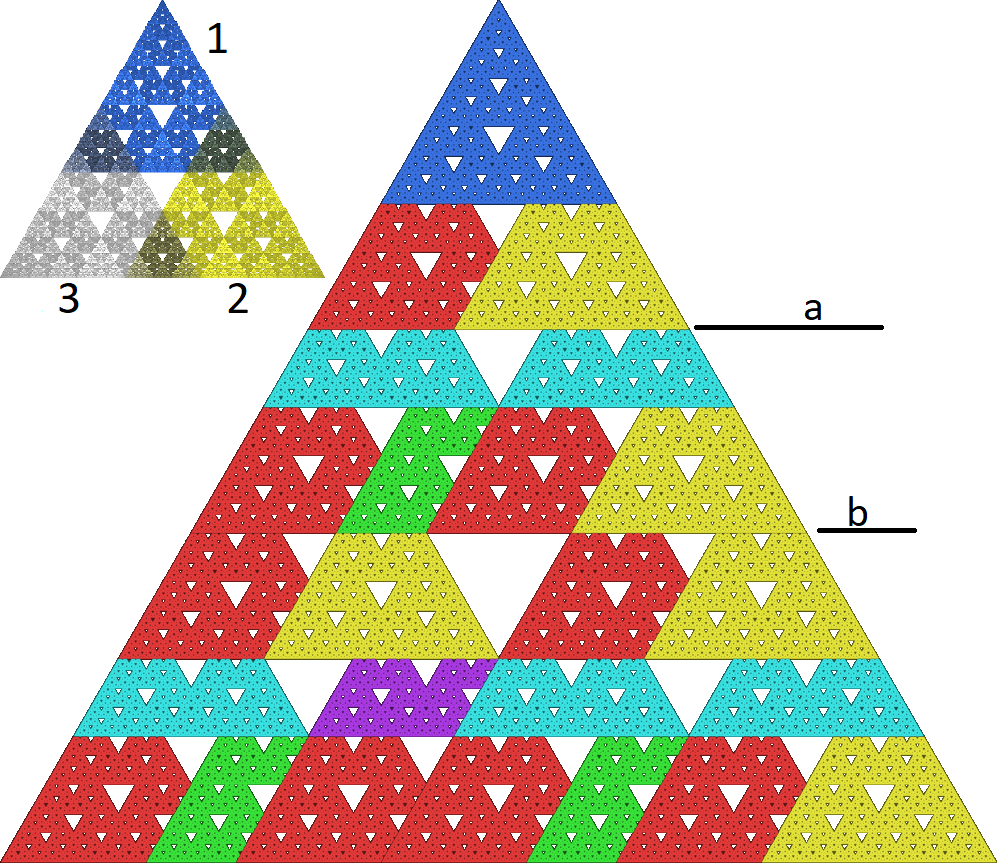} 
\end{center}
\caption{The golden triangle with overlaps, and as a non-overlapping GIFS on level 3.}\label{top2}
\end{figure}

This was essentially a one-dimensional example, with $f_2$ as dummy variable for better visualization. The proper golden triangle in Figure \ref{top2} was introduced by Sidorov et al. \cite{BMS} where the dimension was calculated with an infinite IFS. In the complex plane, the IFS of the symmetric golden triangle can be written as
\begin{equation}
f_1(z)=tz+tv ,\ f_2(z)=tz+t ,\ f_3(z)=tz \ \mbox{ with } v=\textstyle\frac12 (1+i\sqrt{3}) \mbox{ and } t=\frac12 (\sqrt{5}-1)\ .
\label{goldifs}\end{equation}

In establishing the GIFS we follow the principle that each point should be represented by its smallest address with respect to the lexicographic order of the pieces \cite{BB23}.  Thus $B_1=A,$ but $B_2$ must be reduced so that $f_1(B_1)$ and $f_2(B_2)$ do not overlap. From $B_3$ we have to subtract two overlaps corresponding to two neighbors with smaller index. 

We now explain the GIFS. The construction is derived in the next section.
Colors and notation of Figure \ref{top2} follow the first example. The uppermost three triangles above line $a$ illustrate the equation  \[ B_1=f_1(B_1)\cup f_2(B_2)\cup f_3(B_3), \] where the yellow set $B_2$ is pruned at the upper part and the red set $B_3$ at the upper and right parts. The second level subdivision is above line $b.$ Between lines $a$ and $b$ it can be seen how $B_2$ and $B_3$ are subdivided with the help of new sets $B_4,B_5$ drawn in turquoise and green, respectively:
\[ B_2=f_1(B_4)\cup f_2(B_2)\cup f_3(B_3)\, , \ B_3=f_1(B_4)\cup f_2(B_5)\cup f_3(B_3)\, . \]
In the third level we see how $B_4$ and $B_5$ are decomposed. For $B_5$ we need the purple set $B_6.$ The subdivision of $B_6$ would show up on level 4. Here it can be seen at the triangle's left corner.
\[ B_4=f_2(B_2)\cup f_3(B_3)\, , \ B_5=f_1(B_6)\cup f_3(B_3)\, , \ B_6= f_2(B_5)\cup f_3(B_3) \]
These sets have only two pieces since a whole triangle piece was cut off. The system of six set equations together with the IFS \eqref{goldifs} determine the GIFS which clearly fulfils the open set condition \eqref{gosc}. The open sets $U_k, k=1,...,6$ can be taken as interiors of the convex hulls of the respective sets $B_k.$

To calculate the dimension, consider the 0-1-matrix $M=(m_{k\ell})$ with $m_{k\ell}=1$ if $B_\ell$ appears on the right side of the equation for $B_k.$ All attractors $B_k$ have the same dimension $\beta$ since $M$ is irreducible \cite{MW}. The characteristic polynomial of $M$ is $p(\lambda)= \lambda(\lambda-1)(\lambda+1)(\lambda^3-3\lambda^2+3),$ with spectral radius $\lambda_{\rm max}\approx 2.532.$ So $\beta=\log\lambda_{\rm max}/-\log t \approx 1.9306.$ The trick with $\mu (f_i(A))=t^\beta=s$ in the example above gives $1=3s-3s^3$ since $f_i(A)\cap f_j(A)$ is a piece of level 3. This implies $s=1/\lambda_{\rm max}.$ In \cite{BMS}, an infinite IFS leads to the polynomial, and it is mentioned that $s=(2/\sqrt{3})\cos (7\pi/18) .$

\section{The overlap graph} \label{ft}
Let an IFS $F=\{ f_1,f_2,...,f_m \}$ with overlapping attractor $A$ be given. The numbering of the maps can be chosen arbitrarily. It will be important since for overlapping  pieces $f_i(A), f_j(A)$, we shall cut away the intersection from the piece with larger index in the lexicographic order.

Around 2000, a data structure was developed which controls the intersections of pieces $f(A)\cap g(A)$ with $f,g\in F^*$ \cite{A02,Ba97,DKV,Lal97,LNR,NW,ScT,SW}.  For many examples, the standardized intersections $f^{-1}g(A)\cap A,$ the outer boundary sets, form a graph IFS.  Moreover, it is easier to calculate with mappings $f^{-1}g$ than with boundary sets. This led to the neighbor graph which has vertices $h=f^{-1}g$ with $f,g\in F^k$ for some $k\in\NN$ for which $f(A)\cap g(A)\not=\emptyset .$  Such maps $h$ are called proper neighbor maps. Each of these maps represents the relative position of two intersecting pieces of $A,$ up to similarity. Since pieces in different levels have different size, the map $f^{-1}$ performs a standardization, so that the first piece becomes $A$ and the second $h(A).$ 

If there exist only finitely many proper neighbor maps, the IFS is said to be of finite type.  The neighbor maps can be defined recursively, starting with the identity map $id$ and repeatedly applying the formula $h'=f_i^{-1}hf_j$  for all $i,j\in \{ 1,...,m\} .$ Whenever this relation holds, we draw an edge with label $(i,j)$ from vertex $h$ to vertex $h'.$ If the IFS is of finite type, the recursive procedure will be completed after finitely many steps. For this construction, we refer to \cite{EFG3,EFG4,EFG1,BM09} and recent related work \cite{HR22,Lo,LZ17,ScT,TZ20}. This paper takes the neighbor graph as given and discusses its application.

We do not need the full neighbor graph.  We are only interested in intersections $f(A)\cap g(A)$ which contain complete overlaps. That is, there exist  $f',g'\in F^*$ such that 
$ff'(A)=gg'(A).$ In this case $f'^{-1}f^{-1}gg'=id$  which implies that in the neighbor graph there is a path of edges from $h=f^{-1}g$ to $id.$  Thus we consider the subgraph of the neighbor graph determined by all vertices $h$ which admit a directed path to $id.$ This graph will be called the {\it overlap graph.} 

\begin{Proposition}
The overlap graph is obtained from the neighbor graph by matrix operations.
\label{nbo}\end{Proposition}

{\it Proof. }
Let $M$ be the $n\times n$ adjacency matrix of the neighbor graph: $m_{k\ell}=1$ if there is an edge from vertex $k$ to $\ell ,$ and zero otherwise. Let the initial vertex $id$ correspond to $k=1.$ Determine the matrix $N=M+M^2+...+M^{n-1}.$ Note that $n_{k\ell}$ is the number of directed paths of length $\le n-1$ from vertex $k$ to vertex $\ell .$ The adjacency matrix $O$ of the overlap graph is the submatrix of $M$ given by the rows of $N$ with non-zero first entry, and the corresponding columns. The labels of edges remain unchanged.  

For large $n,$ we avoid overflow and reduce the number of matrix multiplications from $n$ to $2\log_2 n$ as follows. Start with $B=N=M$ and let $k$ be the integer part of $\log_2 n.$  Then repeat $k$ times $ N=\mbox{min}(B*N+N,1)\, ; \  B=\mbox{min}(B*B,1)$ where the minimum is taken in every cell of the matrix.\hfill $\Box$

\begin{figure}[h!t] 
\begin{center}
\includegraphics[width=0.54\textwidth]{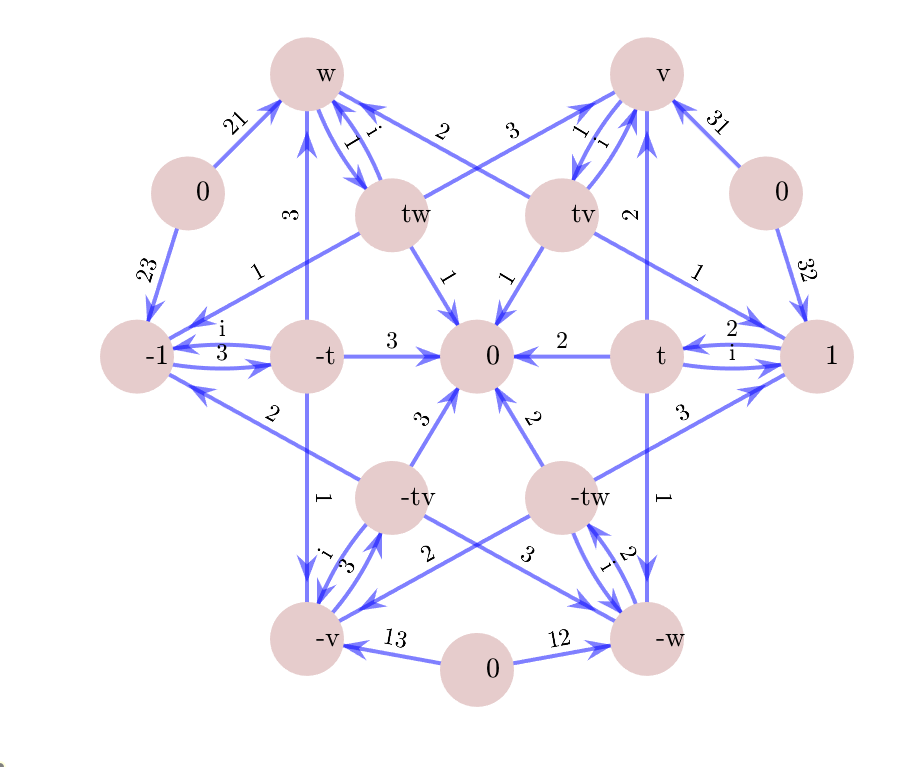}  
\includegraphics[width=0.45\textwidth]{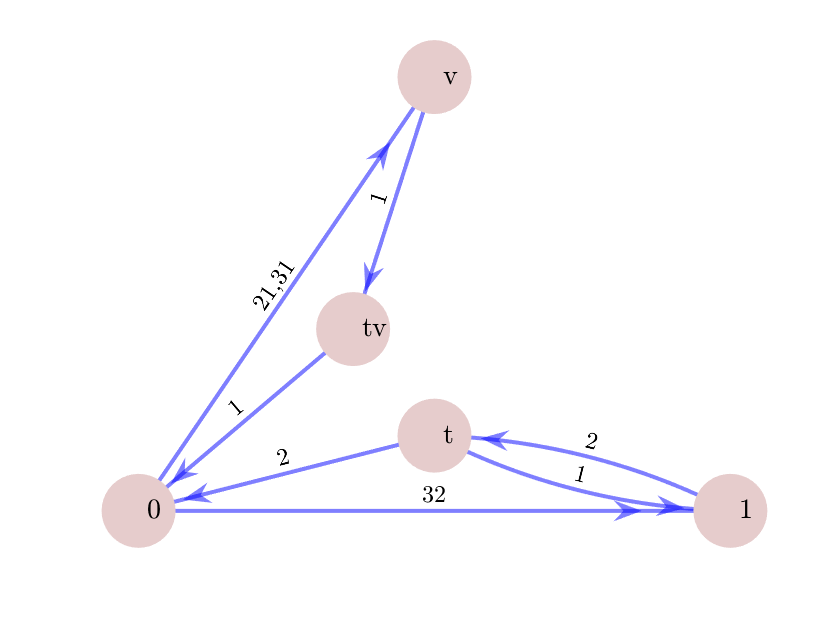} 
\end{center}
\caption{The overlap graph for the golden triangle. To minimize intersections of edges, the initial vertex 0 was drawn three times, and once as terminal vertex in the middle. For the construction of the GIFS, it is sufficient to consider the small subgraph on the right.}\label{top2g}
\end{figure}  

Before we go to details, we briefly review the structure of the overlap graph. Each vertex represents a standardized overlap set $A\cap h(A)$ and the corresponding neighbor map $h.$ There is one initial vertex $id$ or 0, which represents the complete overlap $A=h(A).$ From the initial vertex arise edges with label $(i,j)$ for all $i,j\in\{ 1,...,m\}$ for which $f_i(A)$ and $f_j(A)$ overlap. For another vertex which represents $h$ and the overlap $D=A\cap h(A),$ there arise edges with label $i$ for the $i$ with nonempty overlap $D_i=f_i(A)\cap D .$ The terminal vertex of such an edge is the standardized overlap of the $i$-th piece, $E_i=f_i^{-1}(D_i)$ which is larger than $D_i,$ and $D=\bigcup_i f_i(E_i).$

The overlap graph of the golden triangle is shown in Figure \ref{top2g}. The full neighbor graph has 18 more vertices describing intersections of neighboring pieces in a single point. 
For the IFS \eqref{goldifs}, all neighbor maps are translations. So vertices are denoted by translation vectors. There are translations by the sixth roots of unity, called $\pm 1, \pm v, \pm w$ with $w=v-1,$ which correspond to overlaps in  a small triangle. Translations by the smaller vectors $\pm t, \pm tv, \pm tw$ with  $t=\frac12 (\sqrt{5}-1)$  correspond to large triangular overlaps. It will be sufficient to label edges with a single number $i$ instead of $(i,j).$ Only for the edges starting in $id=0$ we kept the full notation. The letter $i$ itself means ``1 or 2 or 3''. The vertex  $id=0$ was drawn at four places in order to avoid more intersections of edges. 

\begin{figure}[h!t] 
\begin{center}
\includegraphics[width=0.41\textwidth]{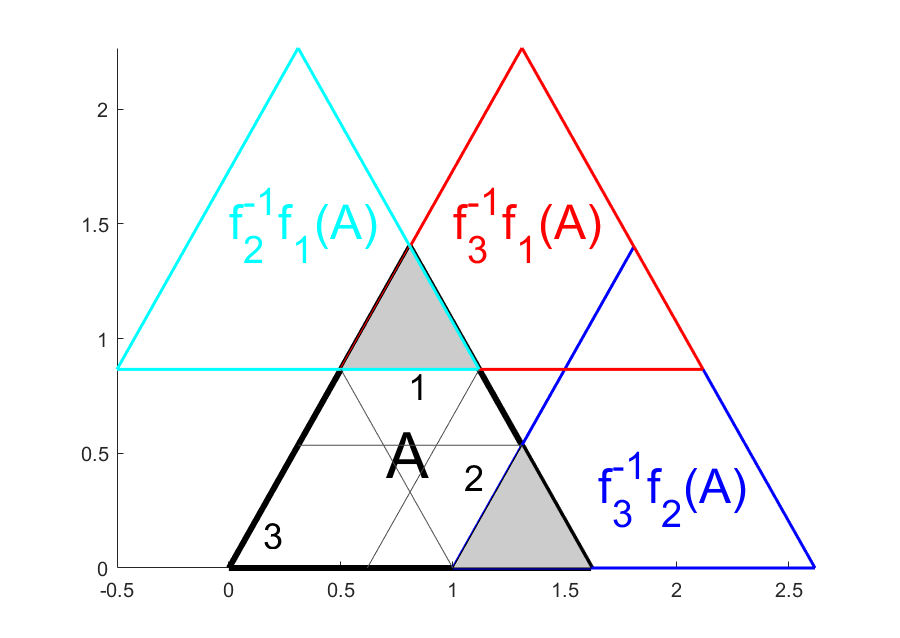}\ 
\includegraphics[width=0.46\textwidth]{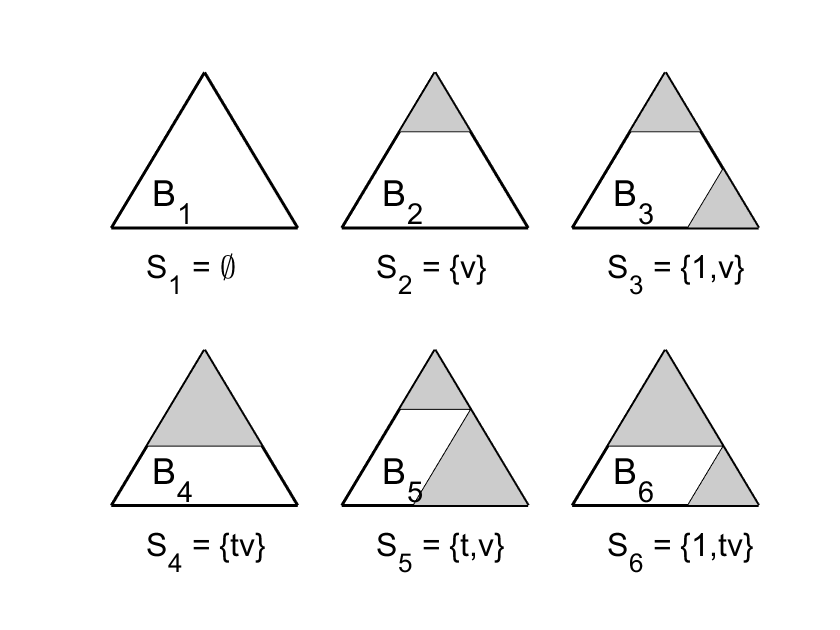}
\end{center}
\caption{Left: The initial overlaps for the golden triangle. Right: The six attractors $B_k.$  Cut overlaps are shaded.}\label{top3}
\end{figure}  

Figure \ref{top3} illustrates the neighbor maps $h=f_i^{-1}f_j$ which correspond to edges starting at the initial vertex 0. We assume $i>j$ since only in this case the neighbor has a lexicographically smaller address and the piece under consideration must be reduced.  
The side length of the triangle $A$ is $t+1=1/t\approx 1.618,$ as indicated in Figure \ref{top3}. From \eqref{goldifs} immediately follows  $f_3^{-1}f_2(z)=z+1$ and $f_3^{-1}f_1(z)=z+v.$ In Figure \ref{top2g} this gives the two edges labelled $32$ and $31$ from 0 on the upper right to the vertices 1 and $v.$ The translation vectors 1 and $v$ represent the two shaded overlap triangles which every piece of type 3 must have with its immediate neighbors. They must be cut away since the neighbors have smaller address. We called the piece $B_3,$ and now we associate it with the overlap set $S_3=\{ 1,v\}$ while the set $B_1=A$ has no overlap so that $S_1=\emptyset .$  For the piece $B_2$ we cut away only the upper small triangle which corresponds to the overlap set $S_2=\{ w\}$ since $f_2^{-1}f_1(z)=z+w.$ 

Because of symmetry, $v$ and $w$ describe the same overlap set. When we put them together, we can use the small subgraph in Figure \ref{top2g} for the following discussion.
We interpret $1$ and $v$ as `small overlaps at the right and the top', shaded in Figure \ref{top3}, and $t,tv$ as `large overlaps at the right and the top', respectively. 

The attractors $B_1, B_2,$ and $B_3$ were obtained from initial edges of the overlap graph. The corresponding sets $S_k$ contain the terminal vertices of these edges. Each vertex represents a neighbor map $h,$ in our case a translation by $1,v,$  or $w,$ which is applied to cut away an overlap. Formally
\[ B_k={\rm cl }\{ A\setminus \bigcup_{h\in S_k}h(A)\}\, , \]
where cl denotes topological closure.  This is the basic formula for all $B_k.$ 
Non-initial edges determine the overlaps of subpieces from the overlap of a piece. As an example, we establish the equation for $B_3$ with $S_3=\{ 1,v \} .$ For piece 1 of $B_3,$ we study edges with label 1 and initial vertex in $S_3.$ There is only one: from $v$ to vertex $tv.$ The small overlap $v$ at the top of $B_3$ becomes a large overlap $tv$ of its piece 1. Thus piece 1 of $B_3$ is determined by the new vertex set $S_4=\{ tv\} ,$  corresponding to $B_4$ in Figure \ref{top3} and to the turquoise set in Figure \ref{top2}.

For piece 2 of $B_3$ we consider edges with label 2 from $S_3.$ The only edge from 1 to $t$ indicates that the small overlap at the right of $B_3$ becomes a large overlap of the piece 2. Moreover, piece 2 has the overlap $v=w$ at the top, since it must respect piece 1 inside $B_3.$ Thus piece 2 has the new vertex set $S_5=\{ t,v\}$ corresponding to the
green set $B_5$ in Figure \ref{top2}.  Since there are no edges with label 3 from $S_3,$ piece 3 of $B_3$ has only the overlaps $S_3=\{ 1,v\}$ imposed by pieces 1 and 2. So piece 3 corresponds to $B_3.$ The resulting equation is $B_3=f_1(B_4)\cup f_2(B_5)\cup f_3(B_3).$ 

Now consider $B_4$ with $S_4=\{ tv\}.$ There is an edge labelled 1 from $tv$ to 0. This means that piece 1 of $B_4$ does not exist - it is contained in the overlap which was cut off. Piece 2 of $B_4$ has only the overlap $S_2=\{ v\}$ (caused not by piece 1, but by the cutoff). Piece 3 again corresponds to $S_3.$ The equation is $B_4=f_2(B_2)\cup f_3(B_3).$ For $B_5$ with $S_5=\{ t,v\}$ piece 2 disappears, and piece 1 has overlaps $S_6=\{ 1,tv\}$ corresponding to another attractor $B_6.$ Overlap 1 at the right of piece 1 is caused by cutting off piece 2. In this way we find the system of equations directly from the overlap graph.

\section{Cut overlaps and neighborhoods} \label{comp}
\begin{figure}[h!t] 
\begin{center}
\includegraphics[width=0.99\textwidth]{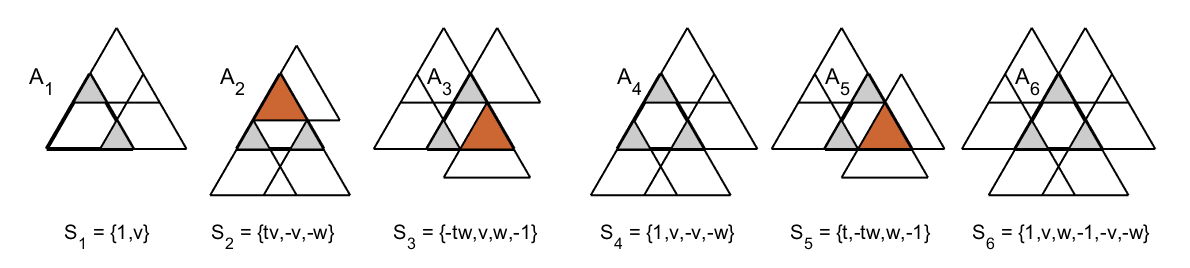}
\end{center}
\caption{The six nontrivial neighborhood types of pieces in the golden triangle. Reflected and rotated neighborhoods have the same type. The brown pieces represent complete overlaps and cause rational entries in the substitution matrix.}\label{top4}
\end{figure}

Our approach is quite similar to the work of Ngai and Wang \cite{NW} which bases on the idea of neighborhood graph developed by Lalley \cite{Lal97}.  In both cases, the number of different overlaps must be finite, up to similarity. The overlapping attractor $A$ is replaced by modifications which we may call types.  In our case the types $B_k$ are copies of $A$ with cut overlaps, as can be seen in Figure \ref{top3}. In \cite{Lal97,NW} neighborhood types $A_k$ are considered. These are pieces of $A$ with a certain environment inside $A,$ that is, neighbors overlapping them in a specific way.  In both approaches, types are determined by sets of overlaps. For the golden triangle, a piece can have two, three, four, or six neighbors.  In the case of four neighbors, three different configurations are possible. Figure \ref{top4} shows the resulting types numbered 1 to 6. Type 0, a set without neighbors which would correspond to our $B_1,$ was omitted.

Rotated and reflected neighborhood types are considered to be identical. Otherwise we would have 21 types for the golden triangle. Using symmetry, the number of types is reduced. As a consequence, a type is given by a whole symmetry class of overlap sets. $S_k$ in Figure \ref{top4} is only one representative.  In contrast, our approach is symmetry-breaking. The sets $B_k$ do not show much symmetry. Their type is given by their very shape which makes our approach simple.

For each type $A_k$ or $B_k,$ the types of the pieces $f_i(A_k)$ or $f_i(B_k)$ are determined for $i=1,...,m$ in both approaches. In this way, Lalley, Ngai and Wang obtain a neighborhood graph while we obtain a system of GIFS equations. In both cases, this results in a ``substitution matrix''. The Hausdorff dimension of $A$ is $\beta=\log \sigma /-\log r\ $ where $\sigma$ denotes the spectral radius of the substitution matrix. The GIFS structure is more powerful since it determines various other topological and measure-theoretical properties of $A.$ One can ask whether the neighborhood approach also leads to a GIFS structure.

\begin{table}[h!]
\begin{center} $\left(\begin{array}{cccccc}
1&1&0&0&0&0\\
0&1&1&1/2&0&0\\
0&0&1&1/2&1&0\\
0&2&0&0&1&0\\
0&0&2&0&0&1/3\\
0&0&0&0&3&0\\
\end{array}\right)$ \end{center}
\caption{The rational ``substitution matrix'' for the neighborhood types $i=1,...,6$ of the golden triangle. Row $i$ counts the successors of a specimen of type $i$ according to successor type.}\label{ratmx}
\end{table}

In general, this is not the case.  The problem is that a substitution matrix with integer entries need not exist, not even for the golden triangle. Apparently, this has gone unnoticed until now. The problem comes with the neighborhoods of pieces which represent complete overlaps, shaded brown in Figure \ref{top4}. This piece in $A_2$ is a successor of type $A_4.$ And it is also a successor $A_4$ of the upper neighbor of the parent. Each $A_4$ piece either has two $A_2$ parents, or two $A_3$ parents, and no other predecessors. Thus each $A_2$ or $A_3$ piece has $1/2$ successor of type $A_4.$

It is possible to circumvent the problem by counting $A_2$ and $A_3$ together, as has been done for other examples in \cite{NW}.  However, type $A_6$ with six neighbors is the child of three different parents $A_5$ and has no other predecessors. There is nothing to count together. The number of type $A_6$ successors of $A_5$ must be 1/3 since three specimen of $A_5$ share one successor $A_6.$ We do not go into details and present only the rational substitution matrix. Its characteristic polynomial is identical to the polynomial $p(\lambda)$ of our substitution matrix at the end of Section \ref{basic}. Apparently, the Hausdorff dimension calculation in \cite{NW} remains true when the substitution matrix contains rational entries.  Further structural information seems out of reach with the neighborhood approach, however.

\section{Combinatorial construction of the GIFS} \label{const}
Now we derive the GIFS equations in the general case. The method, based on calculations with graphs and matrices instead of overlaps and geometry, can be implemented on a computer.  In $\RR^d,$ we are given an IFS $F=\{ f_1,f_2,...,f_m \}$ of similitudes with equal contraction factor and with overlapping attractor 
\begin{equation} 
A=\bigcup_{j=1}^m f_j(A) \ .
\label{AS}\end{equation} 
The indexing of the maps can be chosen arbitrarily but must be fixed for the construction. 
We assume that the overlap graph with vertex set $S$ is given. It can be determined from the neighbor graph by Proposition \ref{nbo}. The construction of the neighbor graph is known  \cite{EFG3,EFG4,EFG1,BM09,HR22,Lo,LZ17,ScT,TZ20}. 

{\it Side remark. }
The overlap graph could be directly constructed from inverse iteration starting in $id.$ Thus we really need only the finite overlap type condition and not the finite type condition which is formally stronger and needed for the neighbor graph. All known examples with overlaps and finite overlap type are finite type. See \cite[Section 2]{EFG3} for a discussion of these subtleties. \vspace{2ex}

Step by step, the equation for $A$ is replaced by a system of equations \eqref{gifs} for the attractors $B_1,...,B_n.$ Each $B_k$ is represented by a set of vertices
$S_k\subseteq S.$  All calculations are done with the $S_k.$ The algorithm is based on \eqref{S1} and \eqref{S2} below. Nevertheless, we add the interpretation for the sets $B_k$ so that the operations for the $S_k$ can be understood. We choose the $S_k$ so that the 
\begin{equation}
B_k={\rm cl }\{ A\setminus \bigcup_{h\in S_k}h(A)\} \, ,
\label{BS}\end{equation}
with cl denoting topological closure, fulfil the GIFS equations in a non-overlapping way. The algorithm consists of an initialization and a recursion. The number of equations is not known at the beginning. The edges of the overlap graph starting in $id$ are labelled $(i,j).$ For all other edges only the first label $i$ is needed.  We use $t\genfrac{}{}{0pt}{1}{i}{\longrightarrow}s$ as shorthand for 
``there is an edge from $t$ to $s$ with label $i$''. \vspace{2ex}
 
{\bf Initialization. } Let $S_1=\emptyset .$ For $i=2,...,m$ let
\begin{equation}
\textstyle  S_i=\{ s\in S\, |\, id \genfrac{}{}{0pt}{1}{(i,j)}{\longrightarrow} s \mbox{ for some }j<i \}\, .
\label{S1}\end{equation}
According to \eqref{BS}, this implies $B_1=A$  and $B_i={\rm cl }\{ A\setminus \bigcup_{j<i} f_i^{-1}f_j(A)\}.$  Thus $f_i(B_i)$ is the closure of $f_i(A)\setminus \bigcup_{j<i} f_j(A).$ In the case that $B_i=\emptyset$ we remove $f_i$ from the IFS without changing $A$ and start the construction with the reduced IFS. By \eqref{AS}, the first equation  
\[ B_1= \bigcup_{i=1}^m f_i(B_i) \]
is now fulfilled, and the overlaps have been removed. Put $k=1$ and $n=m.$ At the end of the initialization, we remove double names for $i=2,...,m.$ When $S_i$ coincides with $S_j$ for some $j<i,$ we use the names $S_j$ and $B_j$ for $S_i$ and $B_i.$ Then we reduce the index of all $S_k, B_k$ with $i<k\le n$ by one and replace $n$ by $n-1.$ Now and in the future $n$ and $k$ denote the current number of different attractors and of established equations, respectively.\vspace{2ex} 

{\bf Recursion. } Let $k=k+1,$ which means we establish a new equation. It has the form
\begin{equation} B_k= \bigcup_{i=1}^m f_i(B_{ik})\ .
\label{Bk}\end{equation}
To this end, we define for $i=1,...,m$ the sets
\begin{equation}
\textstyle S_{ik}= S_i \cup \{ s\in S\, |\, t \genfrac{}{}{0pt}{1}{i}{\longrightarrow}s \mbox{ with }t\in S_k\} \ \mbox{ with } S_i \mbox{ from  \eqref{S1}}.
\label{S2}\end{equation}
According to \eqref{BS}, the set $B_{ik}$ is obtained by subtracting from $A$ the overlaps corresponding to the vertices of $S_{ik}.$ The subset $S_i\subseteq S_{ik}$ cares for the internal overlaps of the pieces of $B_k,$ as in the initialization. The second part of $S_{ik}$ guarantees that $f_i(B_{ik})\subseteq B_k.$ The overlaps $t$ of the big set $B_k$ are expressed as overlaps of its pieces. If there is no edge $t \genfrac{}{}{0pt}{1}{i}{\longrightarrow}s$ then the overlap $t$ does not hit piece $i.$ If there is such an edge, however, the overlap $s$ has to be subtracted from $B_{ik}.$ In this way, \eqref{S2} enforces the validity of the equation for $B_k.$ For certain $i,$ the overlaps can cover the set $A.$ This always happens if $id\in S_{ik}.$ In that case \eqref{BS} says $B_{ik}=\emptyset .$ 

After each recursion step, the terms with $B_{ik}=\emptyset$ are omitted from the equation and the other $S_{ik}$ are renamed for $i=1,...,m.$ If $S_{ik}=S_\ell$ for some $\ell\le n,$ we use the names $S_\ell$  and $B_\ell$ to avoid duplicates. Otherwise, we set $n=n+1$ and introduce the new name $S_n$ for $S_{ik}.$ This means that later we have to build another equation for $B_n.$ The renumbering step transforms the equation \eqref{Bk} into the form \eqref{gifs}.
The recursion has to be repeated until $k=n.$ Then each $B_\ell, \ell=1,...,n$ is expressed by a GIFS equation. 

\begin{Proposition}
\begin{enumerate}
\item[(i) ] The recursion procedure finishes after finitely many steps.
\item[(ii)] For an IFS with the given overlap graph, the constructed system of equations \eqref{Bk} has the unique solution given by \eqref{BS}, $k=1,...,n.$ In particular $B_1=A.$
\end{enumerate}
\label{P2}\end{Proposition}

{\it Proof. }{\it (i).} The finite set $S$ has only finitely many subsets which can appear as $S_k.$ In a related context, this argument was mentioned already 1989 by Thurston \cite[Section 9]{T89}. The theoretical bound $n\le 2^{{\rm card\,}S}$ is huge. 
Our experimental work will show that in practice the number of equations is much smaller.\quad {\it (ii).} For an IFS with contractive maps, any GIFS system of equations \eqref{gifs} has a unique solution consisting of compact nonempty sets $B_k$ \cite{MW}. Since the $S_k$ were constructed so that the $B_k$ with \eqref{BS} fulfil the equations, this must be the solution.  \hfill $\Box$ \vspace{2ex}

\section{Proof of the open set condition} \label{prof}
We now prove Theorem \ref{main}.  The GIFS equations \eqref{gifs} and the family of compact non-empty attractors $B_1,..., B_n$ were constructed above. We must define the open sets $U_1,...,U_n$ and show that they fulfil \eqref{gosc}. The proof for attractors with non-empty interior will serve as a blueprint for the general case.\vspace{1ex}

{\it Proof for ${\rm int\,} A\not= \emptyset .$ } In this case the definition is simple:  $U_k={\rm int\,} B_k.$ Thus $U_1={\rm int\,} A\not=\emptyset .$ At the initialization of the algorithm, before renumbering,
\[ U_k={\rm int\,} A \setminus  \bigcup_{j=1}^{k-1} f_k^{-1}f_j(A ) \ \mbox{ for } k=2,...,m  \]
which implies that the $f_k(U_k)=f_k({\rm int\,} A )\setminus  \bigcup_{j=1}^{k-1} f_j(A )
\subset U_1\setminus \bigcup_{j=1}^{k-1} f_j(U_1)$ are disjoint subsets of $U_1.$ Moreover $U_k\not=\emptyset$ since we assumed $B_k\not=\emptyset .$ So the OSC is fulfilled for $k=1.$  

Any self-similar set $A$ with non-empty interior is regular-closed. That is, ${\rm cl\,}({\rm int\,} A)=A.$ The definition \eqref{BS} of $B_k$ can be reformulated as $B_k={\rm cl }\{ {\rm int\,} A\setminus \bigcup_{h\in S_k}h(A)\}$ because the finite union is compact. Thus all $B_k$ and $B_{ik}$ are regular-closed. 
In an equation \eqref{Bk} for $k>1,$ the set $U_{ik}={\rm int\,} B_{ik}$ can only be empty if $B_{ik}$ is empty. As above, the first part of the definition \eqref{S2} implies that the sets $f_i(U_{ik}), i=1,...,m$ are disjoint.
The second part of \eqref{S2} together with \eqref{BS} guarantees that $f_i(B_{ik})\subset B_k$ and  consequently $f_i(U_{ik})\subset U_k.$ This proves the OSC for the equations \eqref{Bk}. Renaming the $B_{ik}, U_{ik}$ brings them to the form \eqref{gifs} and preserves the OSC. \hfill $\Box$\vspace{2ex}

{\it Proof of the general case. } We work with neighbor maps $h=f^{-1}g$ where $f,g\in F^*.$ Let $N$ denote the infinite set of all such neighbor maps for which $h(A)$ does not contain a whole piece of $A.$ In particular, the (infinite) neighbor graph must not contain a (finite) directed path of edges from $h$ to $id.$ Our open sets should not intersect these non-overlapping neighbors. We use central open sets \cite{BHR,BB21} and define
\begin{eqnarray}
U_1&=&\{ x\in\RR^d\, |\, d(x,A) < \inf_{h\in N} d(x,h(A)) = d(x,{\textstyle \bigcup}_{h\in N} h(A)) \}\ ,\notag \\
U_k&=&\{ x\in\RR^d\, |\, d(x,B_k) < \inf_{h\in N\cup S_k} d(x,h(A)) \}\ 
\mbox{ for } k=1,...,n\, ,\label{Ui}
\end{eqnarray}
and similar for the $U_{ik}$ with $U_k,B_k,S_k$ replaced by $U_{ik}, B_{ik}, S_{ik},$ respectively. Here $d(x,A)=\inf \{ |x-y| \, :\, y\in A\}$ denotes the distance from a set. Proposition \ref{P3} below shows that $U_1$ contains a dense subset of $A.$  Since $B_k\subset A,$ this implies that $U_k$ can only be empty if $B_k\subseteq \bigcup_{h\in S_k} h(A).$ That would mean $B_k=\emptyset .$ This argument also shows $U_{ik}\not=\emptyset$ unless $B_{ik}=\emptyset .$ \vspace{1ex}

For the initialization step we have to show that the $f_k(U_k)$ with $k=1,...,m$ are subsets of $U_1,$ and that they are disjoint. First we note that $U_k\subseteq U_1$  because $B_k\subseteq A$ and so for $x\in U_k$
\[ d(x,A)\le d(x,B_k) < \inf_{h\in N\cup S_k} d(x,h(A))\le \inf_{h\in N} d(x,h(A)) \ .\] Thus it suffices to show that the $f_k(U_1)$ are subsets of $U_1.$ Let $K=\bigcup_{h\in N} h(A).$ Take a point $x\in U_1$ and its image $y=f_k(x).$ Since $f_k$ is a similitude with factor $r,$
\[ d(y,A)\le d(y,f_k(A))=r d(x,A)<r d(x,K)= d(y,f_k(K))\le d(y,K) \ .\]
We used $K\subseteq f_k(K)$ which follows from the fact that each  $h(a)$ with $a\in A$ and $h\in N$ can be written as $f_k(f_k^{-1}hf_j(a'))$ for at least one $j$ for which $a\in f_j(A)$ and $a'=f_j^{-1}(a).$ The neighbor map $h'=f_k^{-1}hf_j$ is a successor of $h$ in the neighbor graph. So it belongs to $N$ whenever $h$ does. We have proved $f_k(U_1)\subseteq U_1.$

Next we show that $f_j(U_j)\cap f_i(U_i)=\emptyset$ for $1\le j<i\le m.$ Due to the definition $B_i={\rm cl }\{ A\setminus \bigcup_{h\in S_i}h(A)\}$ the sets $B_j,B_i$ are closed subsets of $A,$ and the intersection $f_j(B_j)\cap f_i(B_i)$ is relatively small: when  
$f_j(A)\cap f_i(A)$ is an overlap, then $f_i^{-1}f_j$ is in $S_k,$ and 
$f_i^{-1}f_j(A)$ is subtracted before taking the closure to obtain $B_i.$ For $x\in U_i$ the definition states $d(x,B_i)<d(x,f_i^{-1}f_j(A)).$ Applying $f_i$ on both sides we get 
\[ d(f_i(x),f_i(B_i))<d(f_i(x),f_j(A))\le d(f_i(x),f_j(B_j))\ .\]
Thus the points of $f_i(U_i)$ are closer to $f_i(B_i)$ than to $f_j(B_j).$

On the other hand, for an arbitrary small piece $g(A)$ with $g\in F^*$ which does not intersect $f_j(A),$ the neighbor map $h=f_j^{-1}g$ belongs to $N$ since $A\cap h(A)=\emptyset .$ There is a countable family $G\subset F^*$ such that the pieces $g(A)$ do not meet $f_j(A)$ and cover a dense subset of $f_i(B_i).$ According to the definition \eqref{Ui} we have for $x\in U_j$
\[ d(x,B_j)<\inf_{h\in N\cup S_j} d(x,h(A)) \le \inf_{g\in G} d(x,f_j^{-1}g(A))= d(x,f_j^{-1}( {\textstyle \bigcup_{g\in G}} g(A)) \le d(x,f_j^{-1}f_i(B_i)) .\]
Applying $f_j$ on both sides, we see that the points of $f_j(U_j)$ have smaller distance to $f_j(B_j)$ than to $f_i(B_i).$ Together, the two statements imply that $f_j(U_j)$ and $f_i(U_i)$ are disjoint. The OSC for the first equation of the GIFS is verified.\vspace{1ex}

For the equations $B_k=\bigcup_{i=1}^m f_i(B_{ik}), k=2,...,n$ the definition \eqref{Ui} and the proof of the disjointness of the $f_j(U_{jk})$ are the same, only with $B_i,U_i,S_i$ replaced by $B_{ik}, U_{ik},S_{ik}.$  The challenge is to prove $f_i(U_{ik})\subseteq U_k$ for $i=1,...,m.$ Consider a point $y\in U_{ik}$ and its image $x=f_i(y)\in f_i(U_{ik}).$ We show condition \eqref{Ui} for $x.$ The definition \eqref{S2} of $S_{ik}$ guarantees $f_i(B_{ik})\subseteq B_k.$ Thus
\[ d(x,B_k)\le d(x,f_i(B_{ik}))=rd(y, B_{ik})<r\inf_{h\in N\cup S_{ik}} d(y,h(A)) =C \] 
according to definition \eqref{Ui} of $U_{ik}.$ Now $C$ is the minimum of two terms:
\[ r\inf_{h\in N} d(y,h(A)) = \inf_{h\in N} d(x,f_ih(A)) \le \inf_{h\in N} d(x,h(A)) \] because of $f_i(K)\subseteq K$ as shown above, and
\[  r\inf_{h\in S_{ik}} d(y,h(A)) \le  r\inf_{h\in S_{ik}^*} d(y,h(A))= r\inf_{h\in S_{k}} d(y,f_ih(A))\le \inf_{h\in S_{k}} d(x,h(A))\ .\]
Here $S_{ik}^*=\{ s\in S\, |\, t \genfrac{}{}{0pt}{1}{i}{\longrightarrow}s \mbox{ with }t\in S_k\}$ is the second part of $S_{ik}$ in \eqref{S2} which by definition of the neighbor graph refers to the overlaps $f_i^{-1}h(A)$ with $h\in S_k.$ Taking the minimum of the two upper bounds we obtain
\[ d(x,B_k)< C\le \inf_{h\in N\cup S_k} d(x,h(A)) \]
which says that $x$ belongs to $U_k$ by \eqref{Ui}.\vspace{1ex}

Last not least we prove that $U_1$ is nonempty. Fix some $\varepsilon >0$ and let $A^\varepsilon = \{ x\, |\, d(x,A)<\varepsilon \} .$  For $\ell=1,2,...$ and $f\in F^\ell$ let
\[ \gamma (f)= {\rm card}\{ g\in F^\ell\,|\, g(A)\cap f(A^\varepsilon)\not=\emptyset \} \] 
denote the number of pieces of the same size near $f(A).$ Note that equal mappings, as well as overlaps, must be counted only once. Now set $\gamma =\sup_{f\in F^*} \gamma (f).$ Zerner proved that, independently of $\varepsilon ,$ the IFS $F$ fulfils the WSC if and only if $\gamma<\infty$ \cite{Zer}, \cite[Section 4.2]{BSS23}.
For Theorem \ref{main} we assumed that $F$ fulfils the WSC. So $\gamma$ is a positive integer, hence $\gamma=\gamma (\overline{f})$ for some $\overline{f}\in F^*.$ Define the open set
\[ V = \bigcup_{g\in F^*}  g\overline{f}(A^{\varepsilon}) \ .\]
This set was introduced by Schief \cite{Sch} in the OSC setting \cite[Section 4.1]{BSS23}. Obviously $V$ contains a dense subset of $A.$

\begin{Proposition} \qquad
$V\cap h(A)=\emptyset$ for $h\in N\ ,$ and \ 
$V\cap A\subset U_1 \ .$
\label{P3}\end{Proposition}

{\it Proof. } Let $h=\tilde{f}^{-1}\tilde{g}$ with $\tilde{f},\tilde{g}\in F^*.$ 
We assume that there is a point
\[ x\in h(A)\cap V = \tilde{f}^{-1}\tilde{g}(A)\cap {\textstyle \bigcup_{g\in F^*\ }}  g\overline{f}(A^{\varepsilon}) \]
and show that $h$ is not in $N.$ 
We fix a $g\in F^*$ such that $\tilde{f}(x)\in \tilde{g}(A)\cap \tilde{f}g\overline{f}(A^{\varepsilon}).$ Suppose that  $\overline{f}\in F^\ell.$ The definition of $\gamma=\gamma (\overline{f})$ says that there exist mappings $g_j\in F^\ell$ with $g_j(A)\cap\overline{f}(A^{\varepsilon})\not=\emptyset$ for $j=1,...,\gamma .$ Thus there are $\gamma+1$ pieces $\tilde{g}(A)$ and  $\tilde{f}gg_j(A)$ which all intersect $\tilde{f}g\overline{f}(A^{\varepsilon}).$ 

To obtain a contradiction to the maximality of $\gamma ,$ we need pieces of equal size.
Let $\tilde{f}gg_j$ have length $L,$ and $\tilde{g}=f_1f_2\cdots f_k$ with $f_i\in F.$ 
We replace $\tilde{g}$ by $\overline{g}$ with length $L.$ If  $k\ge L,$ we define $\overline{g}=f_1f_2\cdots f_L$ for which $\overline{g}(A)\supseteq\tilde{g}(A).$ If $k<L,$ we select further maps $f_{k+1},...,f_L\in F$ such that $\overline{g}(A)$ with $\overline{g}=f_1f_2\cdots f_L$ contains $\tilde{f}(x).$

The contradiction to the maximality of $\gamma$ implies that $\overline{g}=\tilde{f}gg_j$ for some $j$ between 1 and $\gamma .$ Then $\tilde{f}^{-1}\overline{g}= gg_j$ represents a complete overlap which implies that $h=\tilde{f}^{-1}\tilde{g}$ is not in $N.$ This proves the first assertion of the proposition.

The second assertion is a direct consequence. For $x\in V,$ we have 
$d(x,\bigcup_{h\in N} h(A)) \ge d(x,\RR^d\setminus V)>0.$ Thus any $x\in V\cap A$ belongs to $U_1$ by definition \eqref{Ui}.
\hfill $\Box$\vspace{1ex}

Theorem \ref{main} is proved. The rest of the paper is devoted to concrete examples.

\begin{figure}[h!t] 
\begin{center}
\includegraphics[width=0.5\textwidth]{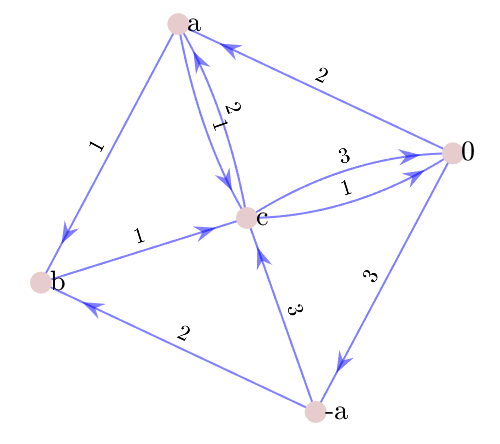}  
\end{center}
\caption{A small overlap graph for which we determine the non-overlapping GIFS.}\label{N24og}
\end{figure}  

\section{A small example} \label{deta}
The GIFS will be constructed for the overlap graph of Figure \ref{N24og} without knowledge of the IFS and the geometry of the attractor. We have three maps $f_1,f_2,f_3.$ There are only two edges starting in the initial state 0, with first labels 2 and 3. Second labels must then be 3 and 2, respectively. Thus only $f_2(A)$ and $f_3(A)$ overlap, and we must subtract the overlap from piece 3 which is lexicographically larger. Thus our first equation is
\[ B_1=f_1(B_1)\cup f_2(B_1)\cup f_3(B_2) \]
where $B_2$ is characterized by $S_2=\{ -a\},$ the endpoint of the initial edge with label 3. For the equation of $B_2,$ we have to consider the successors of $S_2$ in Figure \ref{N24og} with label 1,2, and 3, and only for label 3 we have to add $S_2=\{ -a\}.$
There is no edge with label 1, an edge with label 2 to $b$ and an edge with label 3 to $e.$ According to \eqref{S2} we have $S_{12}=\emptyset =S_1, S_{22}=\{ b\}=S_3$ and $S_{32}=\{ e,-a\}=S_4,$ leading to the equation \[ B_2=f_1(B_1)\cup f_2(B_3)\cup f_3(B_4) \ .\]
From $S_3=\{ b\}$ there is only one edge with label 1 to vertex $e,$ so
$B_3=f_1(B_5)\cup f_2(B_1)\cup f_3(B_2)$
with $S_5=\{ e\}.$ Now consider  $S_4=\{ e,-a\}.$ From $e$ there are edges with labels 1 and 3 to the initial vertex 0. That means that the first and third part of $B_4$ disappear, and only the term with label 2 remains in the equation:
\[ B_4=f_2(B_6) \ \mbox{ with } S_6=\{ b, a\} .\quad \mbox{ Similarly, }
 B_5=f_2(B_7) \ \mbox{ with } S_7=\{ a\} . \]
Using the successors of $S_6$ and $S_7,$ we obtain the equations
\[ B_6=f_1(B_8)\cup f_2(B_1)\cup f_3(B_2) \ \mbox{ and }  \ B_7=f_1(B_8)\cup f_2(B_1)\cup f_3(B_2) \] with $S_8=\{ b,e\},$ and, finally, $B_8=f_2(B_7).$

\begin{figure}[h!] 
\begin{center}
\includegraphics[width=0.4\textwidth]{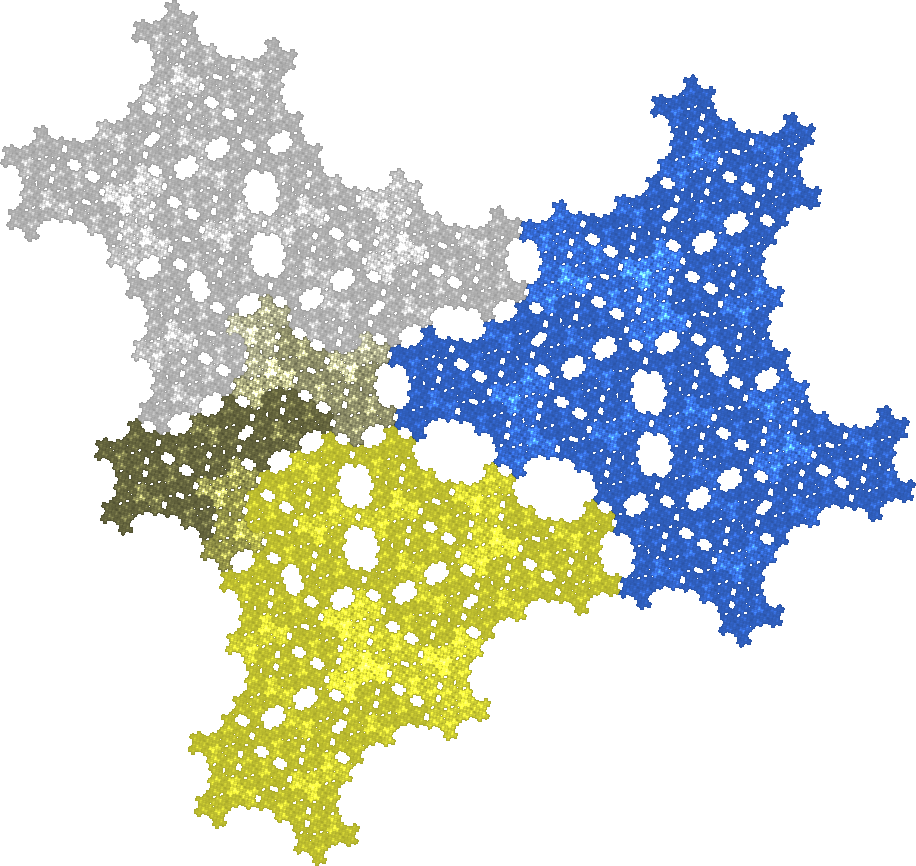}  
\includegraphics[width=0.4\textwidth]{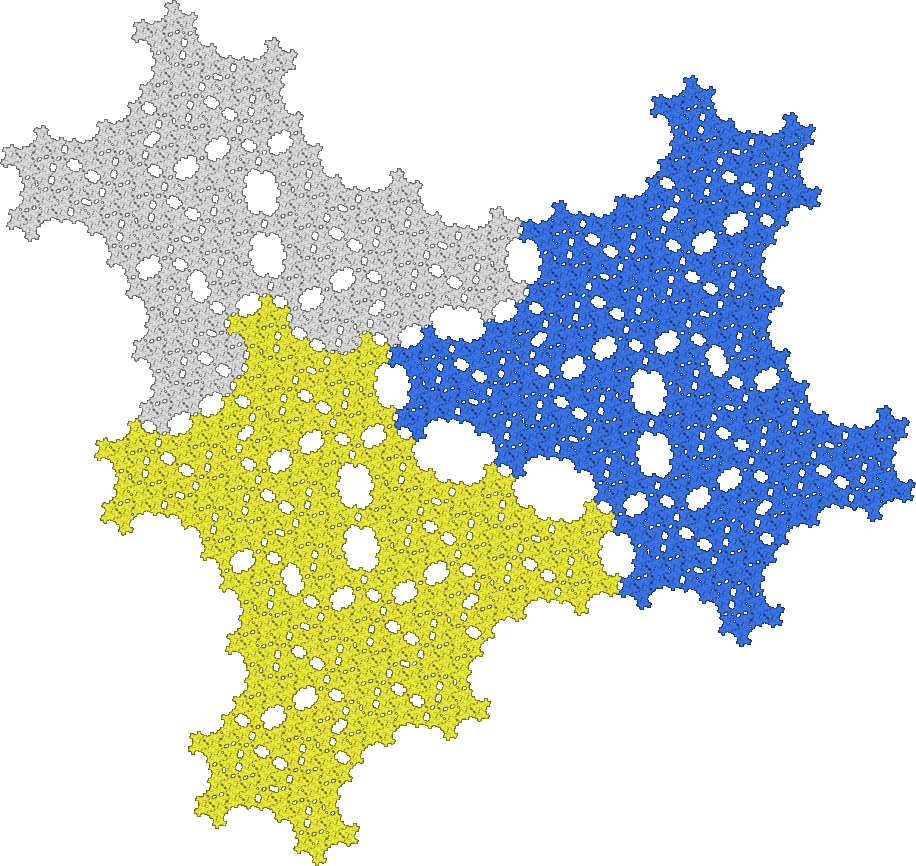} \\
\includegraphics[width=0.4\textwidth]{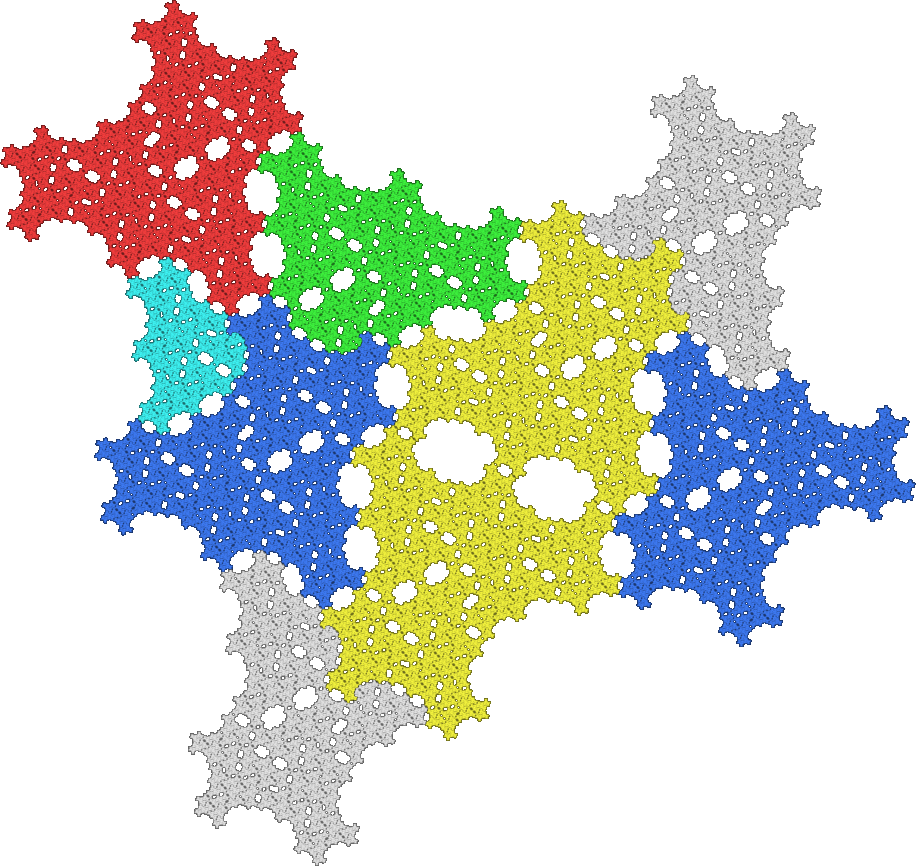} \quad
\includegraphics[width=0.4\textwidth]{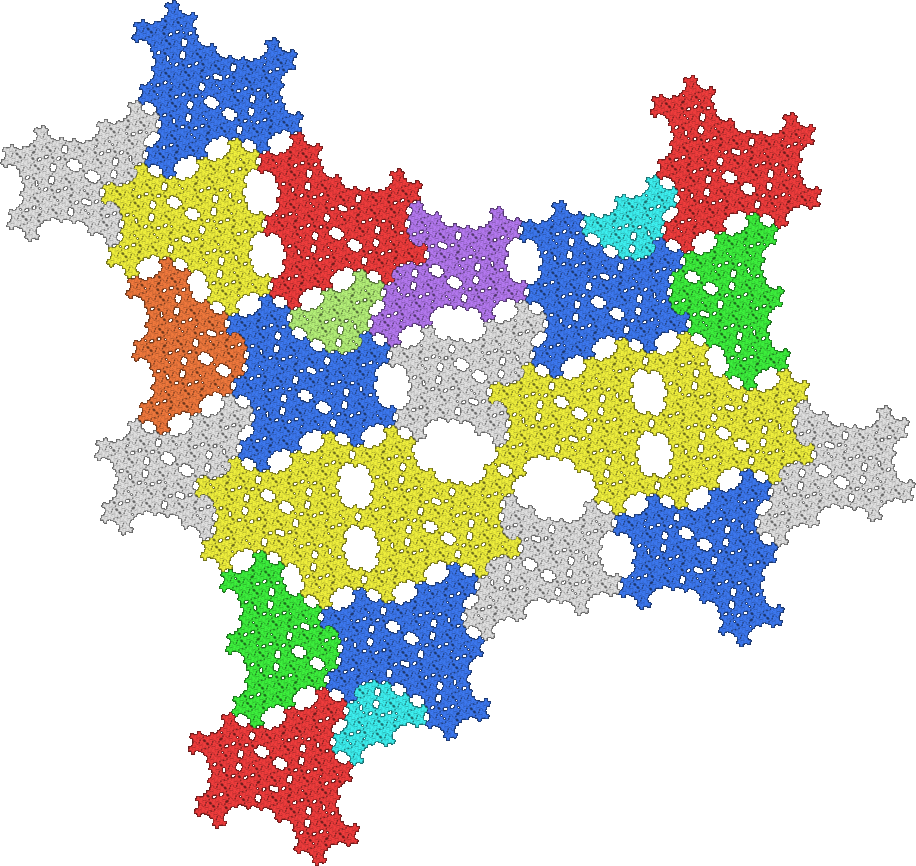}
\end{center}
\caption{Top: A WSC attractor $A$ with cloudy overlaps, introduced in \cite{EFG3}, and the GIFS attractor $B_1$ with clear structure. Bottom: $B_1$ on level 2 and 3. Pieces of type $B_1$ are blue, yellow, and red, $B_2$ gray and pink, $B_3$ green and brown, $B_4$ turquoise and light green.}\label{N24}
\end{figure}  

Since the equations determine their attractors, we have $B_6=B_7$ and $B_8=B_5,$ which then imply $B_5=B_4$ and $B_7=B_3.$ The simplified GIFS system with four attractors consists of the equations for $B_1$ and $B_2$ and
\[ B_3=f_1(B_4)\cup f_2(B_1)\cup f_3(B_2)\ , \quad B_4=f_2(B_3)\ . \]
Figure \ref{N24} shows the overlapping attractor $A$ and the non-overlapping GIFS attractor $B_1$ on levels 1,2, and 3. The coloring of the pieces illustrates the equations for $B_1,..., B_4.$ This example, based on a complex Pisot number of order 4, was introduced and discussed in \cite[Figure 5 and Section 4.2]{EFG3}. 

One can ask why we got a redundant system. The reason is an ambiguity of our coding of overlaps.  Each set $S_k$ characterizes a corresponding overlap set $D_k=\bigcup_{s\in S_k} A\cap s(A)$ where the vertices $s$ of the overlap graph are also neighbor maps. An inclusion $S_k\subseteq S_\ell$  implies $D_k\subseteq D_\ell .$ However, overlap sets $D_k$ can also be subsets of each other, or even coincide, when the corresponding sets of vertices $S_k$ are not at all related.

In the present example, the overlap $D_e=A\cap e(A)$ is very large, containing the first and third piece of $A$ and a part of the second one. It contains the overlaps $D_b, D_a,$ and $D_{-a}.$ Moreover, $D_a=D_b.$ If we knew this, we would have identified $S_4$ with $S_5$, and $S_6$ with $S_3,$ and would have obtained the simple system of four equations immediately.

Actually, we can determine the overlaps $D_s$ for all vertices $s$ mathematically, without further knowledge of the IFS.  The addresses of the points in $D_s$ are given by the sequences of labels of directed edge paths starting in vertex $s$ \cite{EFG4}, using the convention that there are loops from 0 to 0 with label $i$ for $i=1,...,m.$ Instead of infinite edge paths, we can also consider finite edge paths from $s$ to 0, which describe the interior of $D_s$ as a countable union of overlap pieces. This is standard in automata theory, $s$ is the initial  and 0 the final state, and the words labelling the edge paths from $s$ to 0 form a regular language $L_s$ \cite{HU79,epstein}. For $s=e$ we have $L_e=(21)*\{1,3\}$ which denotes the union of pieces 1, 3, 211, 213, 21211, 21213 etc. Since $L_a$ and $L_b$ contain only words starting with 1, it is obvious that they are subsets of $L_e.$ Moreover, $L_b=1L_e$ and $L_a=1L_e\cup 1L_b$ so $L_b\subset L_e$ implies $L_a=L_b.$ It seems interesting that even this simple example leads to questions on regular languages.

\begin{figure}[h!t] 
\begin{center}
\includegraphics[width=0.6\textwidth]{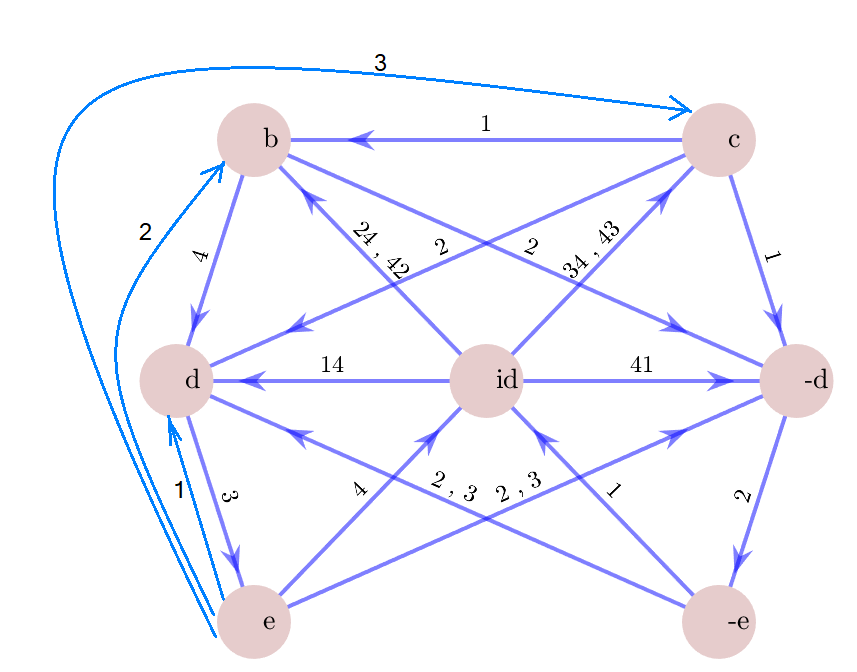}  
\end{center}
\caption{Another overlap graph for which the GIFS can be determined manually. Some edges not needed for the calculation were left, for better visibility.}\label{sq1}
\end{figure}  

\section{Another detailed example} \label{deta2}
The overlap graph in Figure \ref{sq1} describes an example with $m=4$ maps where the overlaps are more separated.   We determine $S_k, k=2,...,m$ according to \eqref{S1}. For $k=2,3$ we find no edges with label 21, 31 or 32. Thus the second and third piece obtain the names $S_1=\emptyset$ and $B_1=A$ after renumbering. The symbols $S_2,B_2$ are needed only for $k=4$ where we have three edges and $S_2=\{ b,c,-d\} :$ 
\[ B_1= f_1(B_1)\cup f_2(B_1)\cup f_3(B_1)\cup f_4(B_2)\, .\]
Next, we determine the equation for $B_2.$ From $S_2,$ outgoing edges with label 1 lead to $b$ and $-d.$ We have a new overlap set $S_3=S_{1,2}=\{ b,-d\}$ and a corresponding attractor $B_3.$ Similarly, outgoing edges from $S_2$ with label 2 lead to $S_4=S_{2,2}=\{ d,-d,-e\} .$ With label 3, however, there are no edges from $S_2,$ so $S_{3,2}=\emptyset =S_1.$ For label 4, we have the overlaps of $S_2=\{ b,c,-d\}$ and $d$ which is approached by an edge from $b.$ With $S_5=\{ b,c,d,-d\} ,$ the equation for $B_2$ is
\[ B_2= f_1(B_3)\cup f_2(B_4)\cup f_3(B_1)\cup f_4(B_5)\, .\]
For $B_3$ with $S_3=\{ b,-d\},$ the pieces 3,4 are the same. But for label 1, no edge can be found, and for label 2, we get the new overlap set $S_6=\{ -d,-e\} :$ 
\[ B_3= f_1(B_1)\cup f_2(B_6)\cup f_3(B_1)\cup f_4(B_5)\, .\]
We continue with the successor vertices of $S_4=\{ d,-d,-e\} .$ For label 1, the successor of $-e$ is $id$ so that the first term of the equation for $B_4$ vanishes. Label 4 leads to $B_2,$ and for labels 2 and 3, we get new sets: $S_7=\{ d,-e\}$ and $S_8=\{ d,e\} :$ 
\[ B_4= \  f_2(B_7)\cup f_3(B_8)\cup f_4(B_2)\, .\]
For $S_5=\{ b,c,d,-d\}$ we get a similar equation as for $S_2=\{ b,c,-d\} .$ The only difference is the edge from $d$ to $e$ with label 3 which leads to $S_9=\{ e\}$ and
\[ B_5= f_1(B_3)\cup f_2(B_4)\cup f_3(B_9)\cup f_4(B_5)\, .\]

Now we can state some properties of the equations. The first term vanishes if $S_k$ contains $-e,$ the last term disappears when $S_k$ contains $e.$ If the last term does not disappear, it will be $f_4(B_5)$ in case $S_k$ contains $b,$ else $f_4(B_2).$ In this way, we derive the equations
\[ B_6= \  f_2(B_7)\cup f_3(B_{10})\cup f_4(B_2)\  \mbox{ with } S_{10}=\{ d\}\ , \ B_7= \  f_2(B_{10})\cup f_3(B_8)\cup f_4(B_2)\ .\]

For $S_8=\{ d,e\}$ we get new sets $S_{1,8}=\{ d\}= S_{10}$ and $S_{3,8}=\{ c, -d, e\}= S_{11}$ and the old set $S_{2,8}=\{ -d, b\}= S_3.$ The last term disappears:
\[ B_8= \  f_1(B_{10})\cup f_2(B_3)\cup f_3(B_{11})\ .\]
Similar results are obtained for $S_9=\{ e\} ,$ only with the new set $S_{12}=\{ c,-d\}$ for $i=3.$ The remaining equations are
\[ B_9= \  f_1(B_{10})\cup f_2(B_3)\cup f_3(B_{12})\ , \ B_{10}= f_1(B_1)\cup f_2(B_1)\cup f_3(B_9)\cup f_4(B_2)\, ,\]
\[ B_{11}= \  f_1(B_{13})\cup f_2(B_{14})\cup f_3(B_{12})\ , \  B_{12}= f_1(B_3)\cup f_2(B_7)\cup f_3(B_1)\cup f_4(B_2)\, ,\]
\[ B_{13}= f_1(B_1)\cup f_2(B_6)\cup f_3(B_9)\cup f_4(B_5)\ ,
\ B_{14}= \ f_2(B_4)\cup f_3(B_8)\cup f_4(B_5)\, .\]\vspace{2ex}

\begin{figure}[h!t] 
\begin{center}
\includegraphics[width=0.4\textwidth]{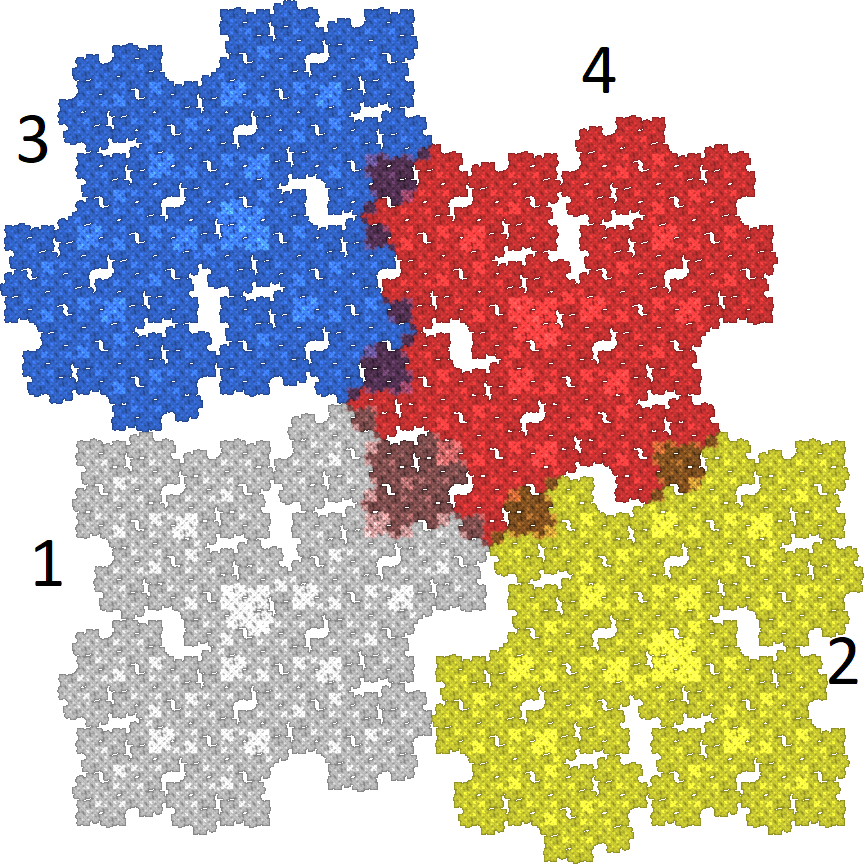}  
\includegraphics[width=0.4\textwidth]{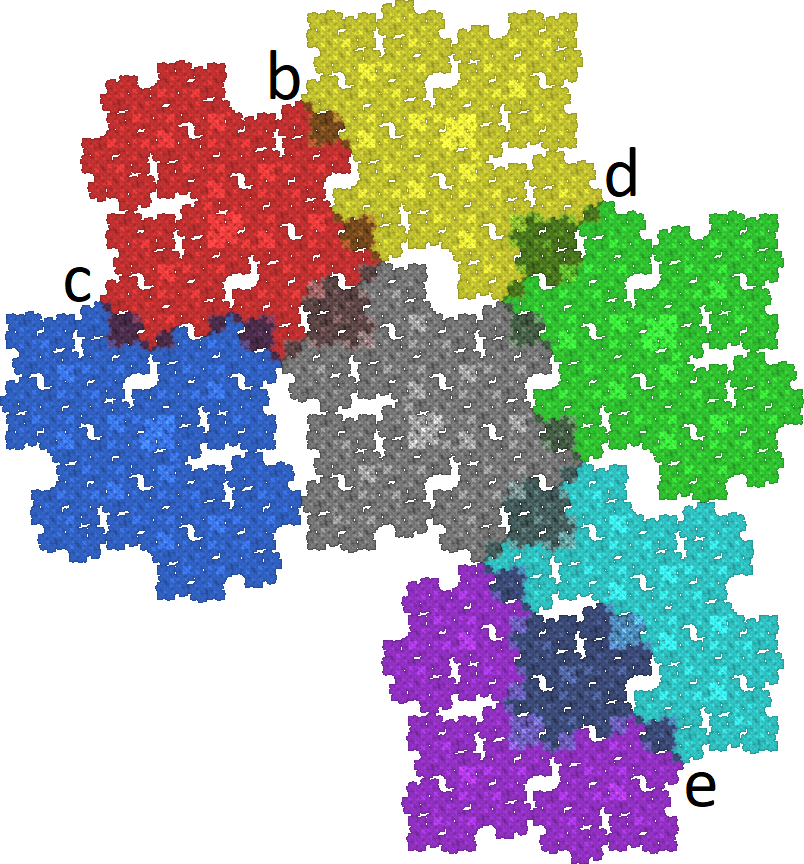} \\
\includegraphics[width=0.4\textwidth]{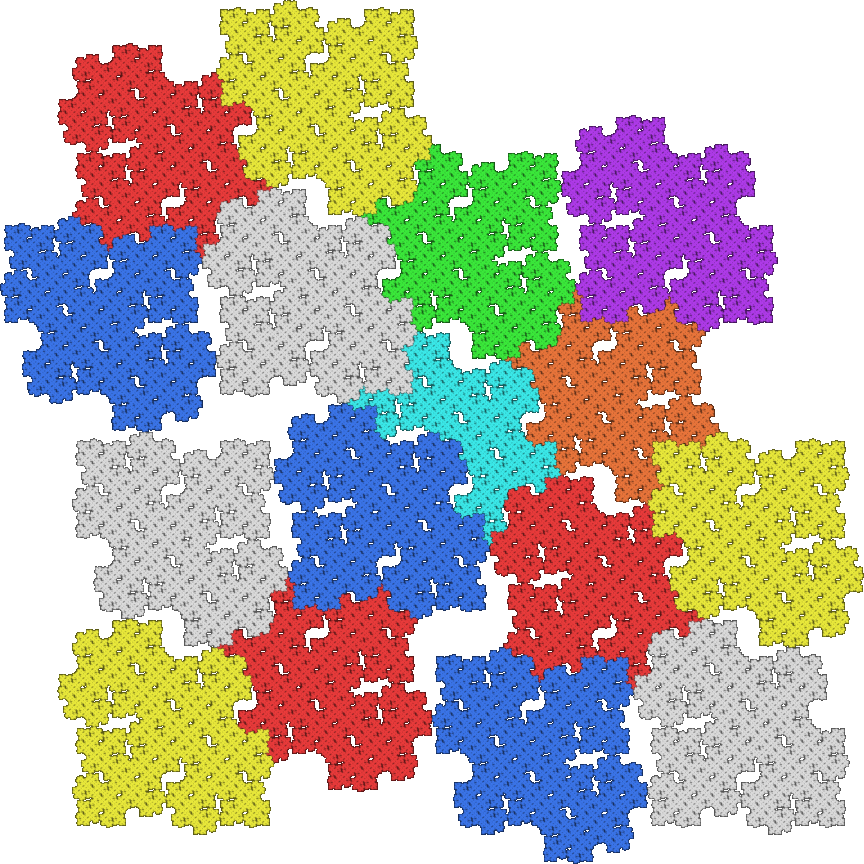} \quad
\includegraphics[width=0.44\textwidth]{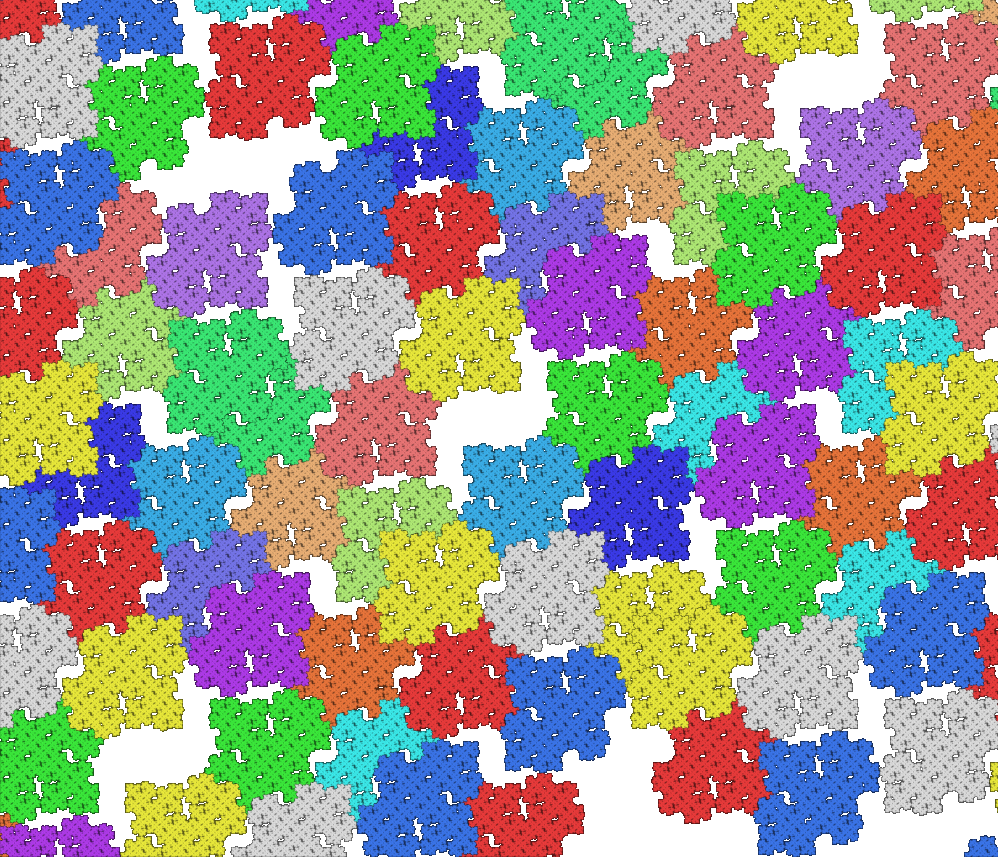}
\end{center}
\caption{Top: the attractor $A$ and its overlaps in the neighborhood of a piece. Bottom: the GIFS attractor $B_1$ on level 2 and a close-up of the fractal tiling generated by the GIFS.}\label{sq2}
\end{figure}  

After establishing the GIFS equations, let us consider the IFS in the complex plane. It is derived from the expanding map $g(z)=2iz,$ with translations and point reflections as `digits':
\[ f_1(z)=\frac{-i}{2}z-1-2i\ , \ f_2(z)=\frac{i}{2}z+1-2i\ , \ f_3(z)=\frac{i}{2}z-2+i\ , \ f_4(z)=\frac{-i}{2}z+1\ . \] 
The IFS fulfils the weak separation condition since the expanding factor $2i$ is a complex Pisot number \cite{EFG3}. The attractor $A$ in  Figure \ref{sq2} can be considered as an overlapping and asymmetric modification of the $2\times 2$ square \cite{EFG1}. The overlaps all occur with piece $A_4,$ while $A_1$ meets $A_2$ as well as $A_3$ in two points only. The overlaps between $A_4$ and $A_2$ as well as $A_3$ correspond to neighbor maps which are point reflections: $b(z)=-z+4, c(z)=-z-2-6i.$ The overlap between $A_1$ and $A_4,$ and between second-level pieces $A_{13}$ and $A_{42},$ is given by the translations $d(z)=z-4+4i, e(z)=z+2+2i,$ respectively. It is correct to write $-d,-e$ for the inverse translations.   Figure \ref{sq2} shows the different types of overlaps in the neighbourhood of a certain piece. The attractor $B_1$ of the constructed GIFS is drawn on level 2, so that the shape of $B_2$ (red), $B_3$ (green), $B_4$ (turquoise), and $B_5$ (green) are visible. The other pieces are copies of $B_1$ (white, yellow, blue, purple).  In the view of the fractal tiling on the bottom right, most of the attractors $B_k$ are represented.

When we tested the method with the reverse ordering of mappings, starting with $S_1=\{ d\}, S_2=\{ b\},$ and $S_3=\{ c\} ,$ we obtained 17 GIFS attractors.  The dimension of the attractors $B_k$ is of course the same, $\beta =\log \lambda/\log 2 \approx 1.9364$ where $\lambda$ is the spectral radius of the $n\times n$ incidence matrix $S$ of the $B_k.$ The matrices have characteristic polynomials with different degree $n,$ but both contain the minimal polynomial $\lambda^5-4\lambda^4-2\lambda^3 +9\lambda^2+4\lambda+2$ of $\lambda$ as a factor.

\section{Computational problems} \label{prob}
\begin{figure}[h!t] 
\begin{center}
\includegraphics[width=0.32\textwidth]{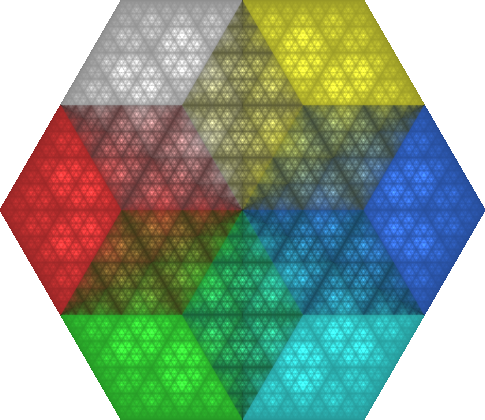}  
\includegraphics[width=0.32\textwidth]{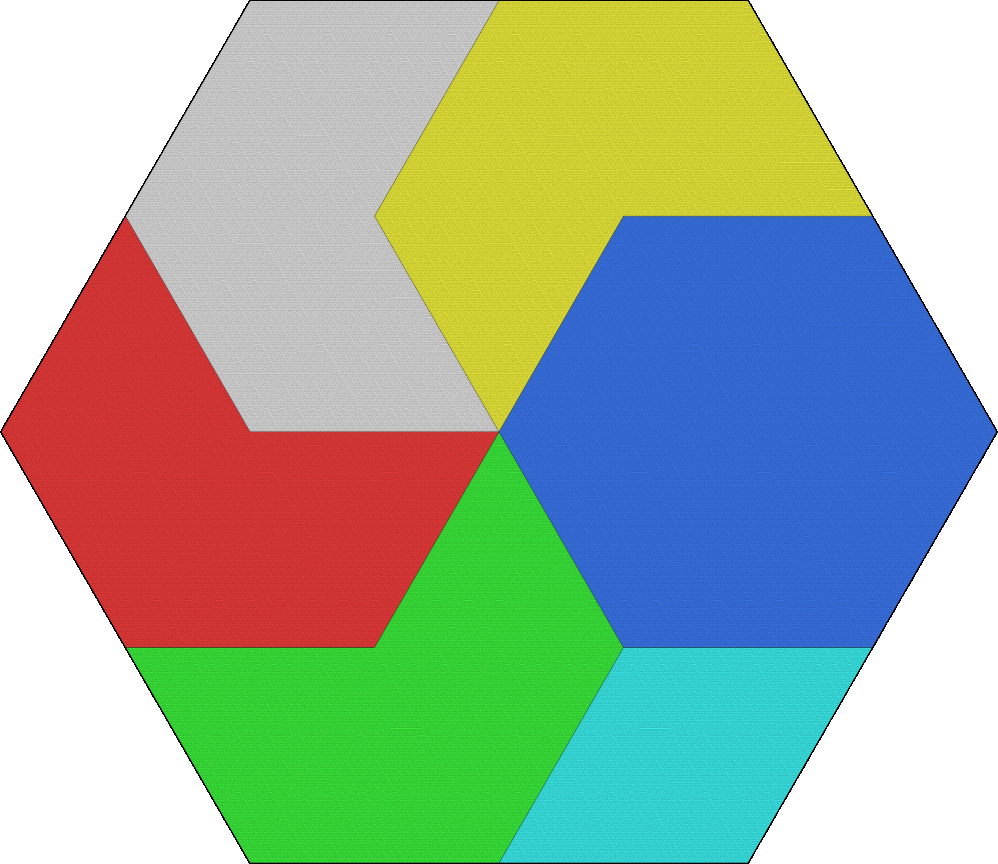} 
\includegraphics[width=0.32\textwidth]{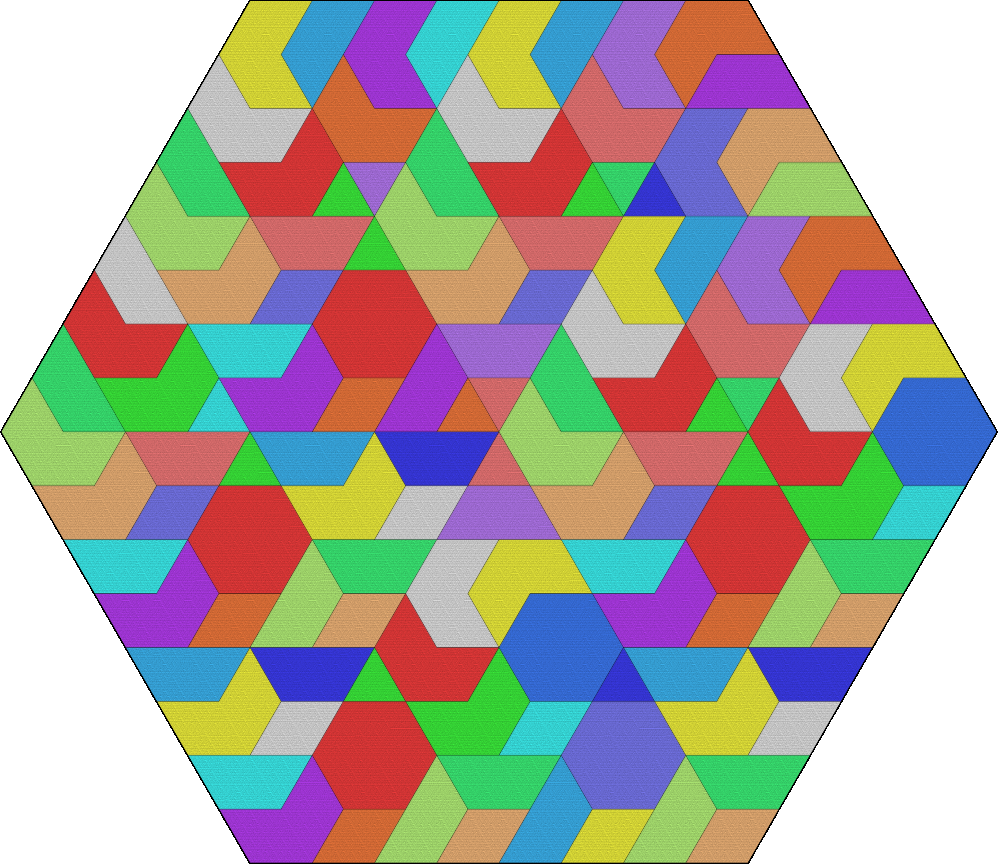} 
\end{center}
\caption{Left: the hexagon as overlapping attractor. The mappings are $f_k(z)= \frac{z}{2} +\omega^{k-1}, k=1,...,6$ where $\omega=\frac12(1+i\sqrt{3})$ fulfils $\omega^6=1.$ Middle and right: the GIFS representation of $B_1$ on level 1 and level 3.}\label{hexa}
\end{figure}

Larger examples are better proceeded by computer.
The GIFS algorithm was implemented in MATLAB and checked for various IFS.   In all cases, the computer-generated GIFS formally fulfilled Theorem \ref{main}.  However, our algorithm does not minimize the number of GIFS attractors, due to the very special form \eqref{gifs} of the equations. In Figure \ref{hexa}, for instance, equivalence of attractors is up to translation. So the four congruent pieces of the picture in the middle are different attractors. Even the translations have a special form so that the two types of triangles visible in the picture on the right represent four attractors. 

It is a mathematical problem to find all possible GIFS representations of some attractor $B_1,$ or at least a minimal one. Our examples indicate that in more complicated examples there are often many GIFS representations. For self-similar sets with the OSC, however, the existence of two truly different IFS representations is rather an exception.  

Our computer experiments also led to computational problems. The algorithm can produce many unnecessary sets $B_k$ and networks of equations for hidden duplicates. It can happen that two vertices of the overlap graph represent the same overlap, as $v$ and $w$ in Figure \ref{top2g} or $a$ and $b$ in Figure \ref{N24og}. A more frequent case is that the overlap of one vertex $s$ is contained in the overlap given by another vertex $t,$ as was explained in Section \ref{deta} in terms of the regular languages $L_s.$ This happened also in Figure \ref{sq1} where $L_e=4\cup 22L_{-e}\cup 32L_{-e}$ and $L_b=43L_e\cup 22L_{-e}$  implies $A\cap b(A)\subset A\cap e(A).$ So $S=\{ e\}$ and $S=\{ e,b\}$ would describe the same attractor. Fortunately, the $S_k$ in this example were not concerned.

In larger examples, the sets $S_k$ may contain more than 10 vertices, and the question then is whether the language $L_s$ for some vertex $s$ is contained in the large language generated by all vertices of $S_k.$ In automata theory, various algorithms have been developed for such problems \cite{BBCF}. The challenge is to give an efficient implementation of a suitable method for our purpose.

In the overlapping hexagon of Figure \ref{hexa} the overlaps are rhombs. Let us number them cyclically 1,...,6. Overlaps 1 and 2 together form a trapezium, half of a hexagon. The $B_k$ generated by the overlaps 1,2,4,5 should then be empty since the overlaps cover the whole hexagon. However, the algorithm does not recognize this fact and establishes an equation. The GIFS solution for this $B_k$ becomes an interval -- the diameter between the two halves -- which then appears as part of other $B_j.$ The overlaps 1,3,5 also cover the hexagon but do not generate such a phantom attractor. Three of 18 attractors were intervals and thus superficial.  

From a theoretical viewpoint, such lower-dimensional attractors are interesting since they are not detected by regular languages. We need the infinite sequences of addresses to verify that the overlaps 1,2,4,5 really cover the hexagon.  However, from a practical viewpoint lower-dimensional attractors are easy to remove. We just take the irreducible part of the incidence matrix $S$ of the $B_k.$ 

For the two examples of the next section, it was enough to adopt a heuristical strategy: take the irreducible part, order the matrix of equations so that identical equations stand next to each other, remove all but one of a group of identical attractors, rename their appearance in other equations, renumber the remaining sets and equations and try to repeat the whole procedure. This was not sufficient for more complicated examples. Even for a WSC Cantor set with only 37 overlaps, after reducing the number of equations from 240 to 185, we are not sure that further reduction is impossible. More systematic computer work is necessary to treat such cases.

\section{Self-replicating tilings} \label{tile}
If the self-similar set $A$ has non-empty interior, it is called a tile. By repeated application of $f^{-1}$ with $f\in F^*$ and fixed point in the interior, we generate a covering of $\RR^d.$ Our GIFS will turn this covering into a proper tiling. Such tilings are known as self-replicating or self-similar tilings \cite{Ke96,T89}, and there is an extensive literature, see for instance \cite{tilingencyc,So20}. Well-known examples are the aperiodic tilings by Robinson, Penrose, and Ammann \cite{Baake2013}, \cite[Chapter 11]{GS}. Other examples, first treated by Rauzy, are obtained from substitutions with a Pisot eigenvalue of the incidence matrix \cite{ABB,AI}. Our algorithm provides an alternative to substitutions and yields new examples. We first construct an overlapping tiling. The finite type property is fulfilled when the expanding factor is complex Pisot and the other data are algebraic integers in the generated number field \cite{EFG3}. We get a covering when the IFS has sufficiently many mappings.  Then we construct the GIFS as above.

\begin{figure}[h!t] 
\begin{center}
\includegraphics[width=0.32\textwidth]{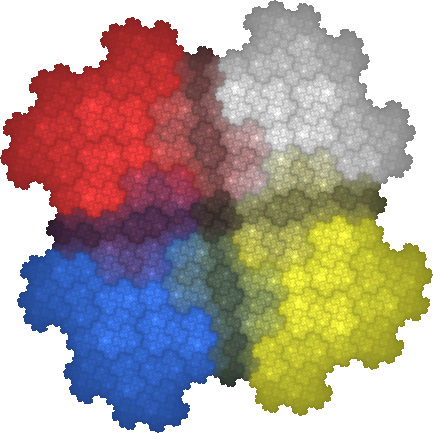}  
\includegraphics[width=0.32\textwidth]{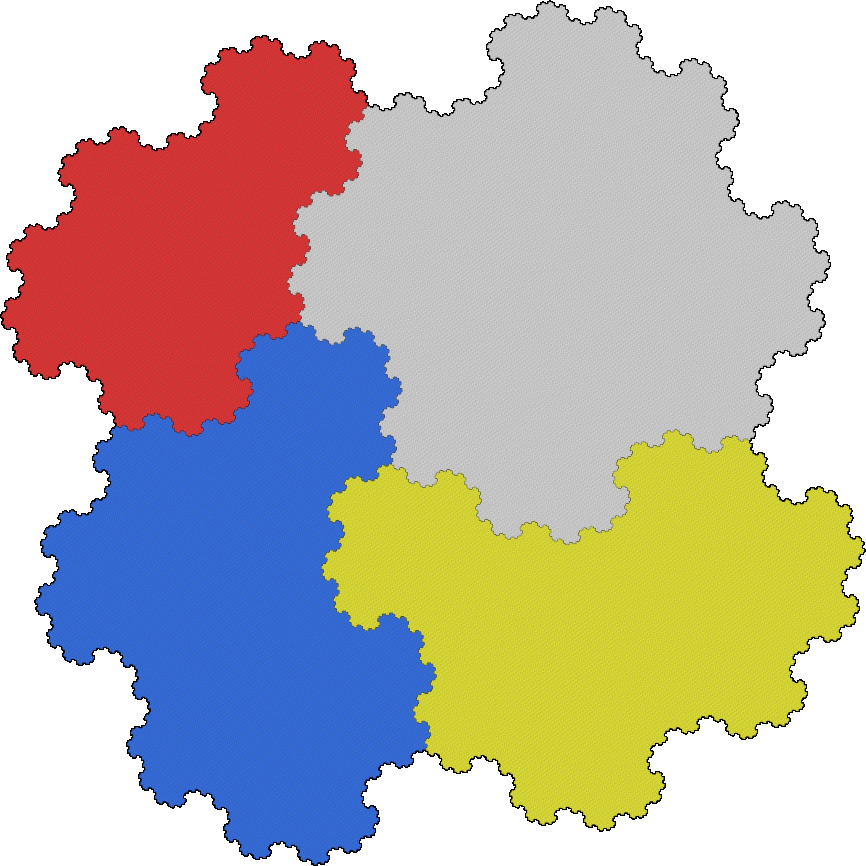} 
\includegraphics[width=0.32\textwidth]{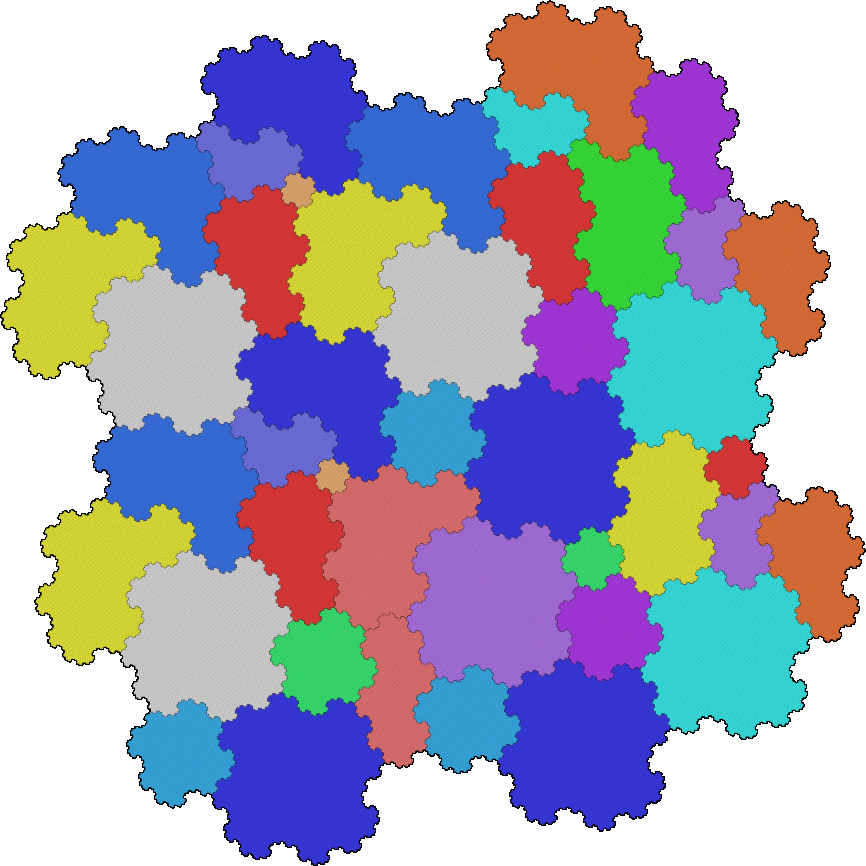} 
\end{center}
\caption{Left: the overlapping tile obtained from a complex Pisot number of degree 4 and translation vectors which form a square.  Middle and right: the GIFS representation of $B_1$ on level 1 and level 3. Thirteen sets $B_k$ are needed.}\label{CP24}
\end{figure}  

For an example, we take the complex Pisot number $\lambda$ with minimal polynomial
$p(\lambda)=\lambda^4+2\lambda^3+4\lambda^2+2\lambda+1.$ Since $p(\lambda)=(\lambda^2+\lambda+1)^2+\lambda^2,$ we have $\lambda+1/\lambda+1=\pm i .$  Thus $i$ is in the generated number field. We choose the IFS  $f_k(z)=\frac{z}{\lambda}+v_k$ so that $v_1=-\lambda, v_2=i\lambda=1+\lambda+\lambda^2, v_3=\lambda, v_4=-i\lambda$ form a square with midpoint zero, generating a  symmetric attractor. Since $|\lambda|\approx 1.7<2,$ the formal similarity dimension is $\log 4/\log |\lambda|>2.6$ and it is quite plausible that the attractor is a tile, even though there are heavy overlaps.

In this example we have 45 proper neighbor maps from which 29 describe overlaps. The theoretical bound for the number of GIFS attractors is of magnitude $2^{28}.$ However, the practical situation is much better. The algorithm determined 51 sets $B_k,$ but the irreducible component of $B_1$ contained only 18 sets. Removing obvious duplicates we arrived at 13 GIFS attractors which are represented in the image on the right of Figure \ref{CP24}. Apparently all $B_k$ are homeomorphic to a disk. They include copies of the symmetric attractor $A$ with four different sizes. Actually, the GIFS system contains three equations with a single term:
\begin{eqnarray*}
B_1= f_1(B_1)\cup f_2(B_2)\cup f_3(B_3)\cup f_4(B_4)\, ,\
B_2= f_1(B_5)\cup f_2(B_6)\cup f_3(B_3)\cup f_4(B_4)\, ,\quad\qquad\\
B_3= f_1(B_1)\cup f_2(B_7)\cup f_3(B_8)\cup f_4(B_4)\, ,\
B_4= f_1(B_1)\cup f_2(B_7)\, ,\ B_6= f_1(B_{10})\cup f_4(B_9)\, , \quad\\
B_5= f_1(B_1)\cup f_2(B_2)\cup f_4(B_9)\, ,\
B_7= f_1(B_{11})\cup f_2(B_6)\cup f_3(B_3)\, ,\ B_8= f_1(B_{12})\cup f_2(B_7)\, , \quad\\ 
B_9= f_2(B_7),\ B_{10}= f_1(B_1)\cup f_2(B_2),\ B_{11}=f_4(B_9),\ 
B_{12}= f_3(B_3)\cup f_4(B_{13}),\ B_{13}=f_1(B_1)\, . 
\end{eqnarray*} 
The main factor of the characteristic polynomial of the incidence matrix of the $B_k$ is $x^8+x^7-3x^6-15x^5-20x^4-15x^3-3x^2+x+1,$ and the spectral radius is a real Pisot number, which within numerical accuracy coincides with $|\lambda|^2,$ confirming the tiling property.

\begin{figure}[h!t] 
\begin{center}
\includegraphics[width=0.32\textwidth]{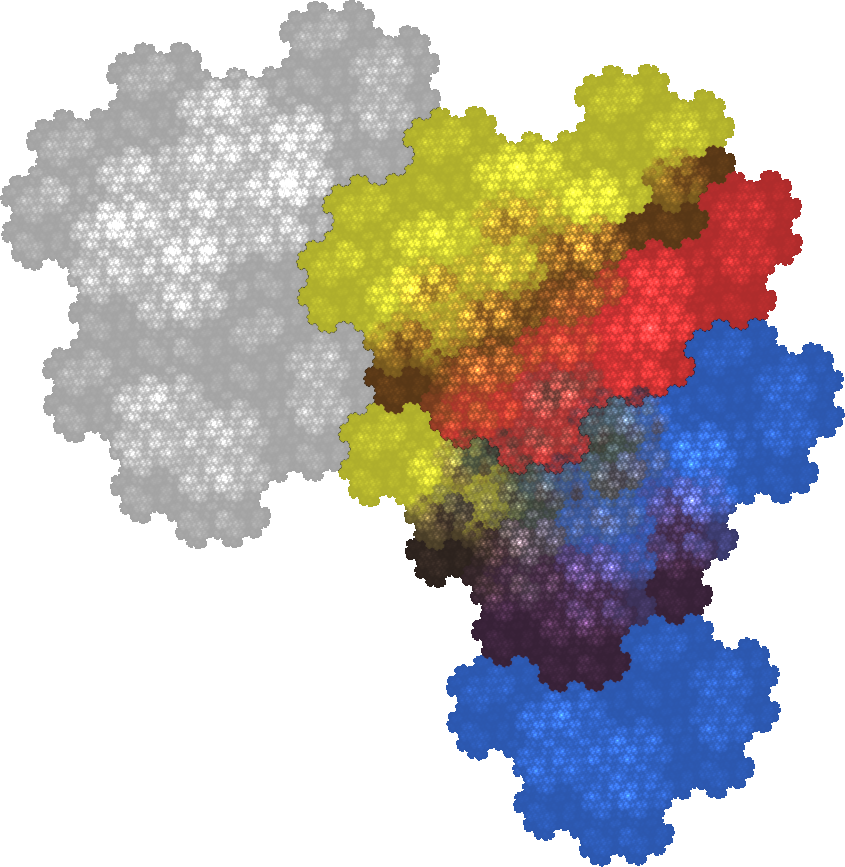}  
\includegraphics[width=0.32\textwidth]{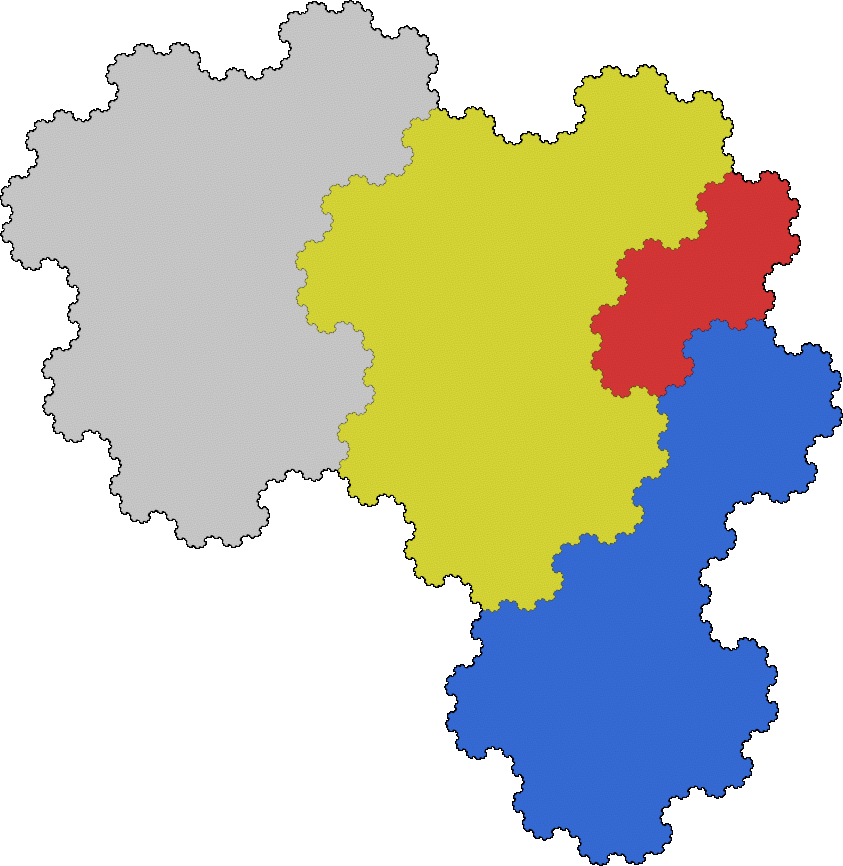} 
\includegraphics[width=0.32\textwidth]{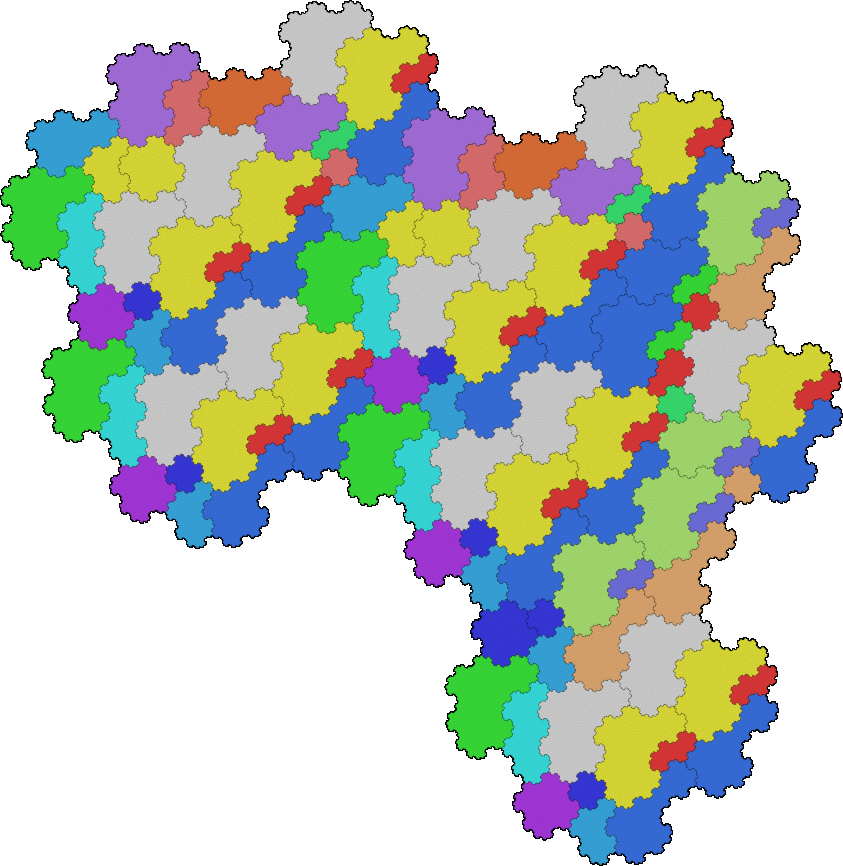} 
\end{center}
\caption{Left: a non-symmetric modification of Figure \ref{CP24} which admits 67 overlaps and still has only 13 attractors $B_k.$  Middle and right: the GIFS representation of $B_1$ on level 1 and 4.}\label{CP24N}
\end{figure}  

Figure \ref{CP24N} shows a non-symmetric example with the same complex Pisot factor. The $v_k$ now are $\lambda , i\lambda -1, \lambda (1-i), \lambda^3+\lambda (1-i).$ There are 97 proper neighbor maps, among them 67 overlaps. Nevertheless it was possible to reduce the initial number of 87 equations to 13 by taking the irreducible part of the matrix and removing overlaps. Among the 13 attractors, three pairs seem to be translations of each other. Such examples give some hope to get small numbers of equations even for complicated overlap structure.

\bibliographystyle{plain} 
\bibliography{lit2}

\end{document}